\documentclass[12pt]{article}
\usepackage[a4paper,left=2cm,right=1.2cm,top=3cm,bottom=4cm]
{geometry}
\usepackage{amsmath,amstext, amsthm, empheq, extpfeil,
amsbsy,amssymb,marvosym,fancyhdr,graphicx,amscd,amsfonts,latexsym,delarray,stackrel,lineno,color,cite,appendix,xcolor,url,hyperref}
\usepackage{wasysym}
\usepackage{subfig}

\usepackage{srcltx}
\usepackage{mathrsfs}
\usepackage{float}
\usepackage[all]{xy}
\usepackage{t1enc}
\usepackage{mathrsfs}
\usepackage{pifont}

\usepackage{mathpple}
\usepackage[T1]{fontenc}

\definecolor{Green}{rgb}{0,1,0}
\definecolor{Blue}{RGB}{0,0,191}
\definecolor{mathmodecolor}{RGB}{0,102,0}
\definecolor{keywordcolor}{RGB}{0,51,151}
\definecolor{sourcebackgroundcolor}{RGB}{255,247,223}
\definecolor{unixagred}{RGB}{255,0,0}
\definecolor{lightgray}{RGB}{191,191,191}
\definecolor{green}{RGB}{1,191,191}

\newcommand*\patchAmsMathEnvironmentForLineno[1]{%
  \expandafter\let\csname old#1\expandafter\endcsname\csname #1\endcsname
  \expandafter\let\csname oldend#1\expandafter\endcsname\csname end#1\endcsname
  \renewenvironment{#1}%
     {\linenomath\csname old#1\endcsname}%
     {\csname oldend#1\endcsname\endlinenomath}}%
\newcommand*\patchBothAmsMathEnvironmentsForLineno[1]{%
  \patchAmsMathEnvironmentForLineno{#1}%
  \patchAmsMathEnvironmentForLineno{#1*}}%
\AtBeginDocument{%
\patchBothAmsMathEnvironmentsForLineno{equation}%
\patchBothAmsMathEnvironmentsForLineno{align}%
\patchBothAmsMathEnvironmentsForLineno{flalign}%
\patchBothAmsMathEnvironmentsForLineno{alignat}%
\patchBothAmsMathEnvironmentsForLineno{gather}%
\patchBothAmsMathEnvironmentsForLineno{multline}%
}

\newtheorem{thm}{Theorem}[section]

\newtheorem{rem}[thm]{Remark}

\def\F{{\mathbb F}}

\def\N{{\mathbb N}}

\def\Q{{\mathbb Q}}

\def\Z{{\mathbb Z}}

\def\aarith{{\mathfrak A}}

\def\cO{{\mathcal O}}

\def\scal{{(\rnt,\cO)}}
\def\dar[#1]{\ar@<2pt>[#1]\ar@<-2pt>[#1]}

\newcommand{\ie}{{\it i.e.\/}\ }

\def\bm2{{\rm Bmod^2}}
\def\b2{{\rm Bmod^{\mathfrak s}}}

\newcommand{\nil}[1]{}

\def\rnt{{[0,\infty)\rtimes{\N^{\times}}}}

\def\aarith{{\mathscr A}}
\def\scal{{(\rnt,\cO)}}
\def\scal1{{\hat \aarith}}
\def\scal2{{\mathscr S}}

\title
{Motivic Rhythms}

\author{Alain Connes}

\begin{document}

\maketitle

\begin{abstract}
In  this article  on mathematics and music, we explain how one can "listen to motives" as rhythmic interpreters.  In the simplest instance which is the one we shall consider, the motive is simply the $H^1$ of the reduction modulo a prime $p$ of an hyperelliptic curve (defined over $\Q$). The corresponding { time onsets} are given by the arguments of the complex eigenvalues of the Frobenius. We find  a surprising relation between mathematical properties of the motives and the ideas on rhythms  developed by the composer Olivier Messiaen.
\end{abstract}
\tableofcontents

\newpage

 The name "motive"  { has a precise  meaning in   mathematics where it was}  coined by Alexander Grothendieck.  { We explain below how  mathematical motives give rise to rhythmic interpreters and illustrate their variety by simple concrete examples}. An ordinary geometric space, such as a compact Riemannian space, has an associated scale in the musical sense as exemplified in the famous lecture ``can one hear the shape of a drum" of Marc Kac { \cite{Kac}}. The association from motives to music explained below is  quite different.  Our thesis is that a motive manifests itself not by a musical scale (which would be a collection of frequencies) but rather by a (periodic) collection of {``time onsets"}. In the simplest instance which is the one we shall consider { as a rich source of concrete examples}, the motive is simply the cohomology $H^1$ of the reduction modulo a prime $p$ of an hyperelliptic curve (defined over $\Q$). The corresponding {``time onsets"} are given by the arguments of the complex eigenvalues of the Frobenius. Our goal which is of an { "illustrative nature"} is to exhibit the simplest motives as rhythmic interpreters. { To an hyperelliptic curve of genus $g$ corresponds for each prime number $p$ not dividing the discriminant a collection of $2g$ time onsets with profound mathematical meaning which repeat in a periodic manner with a frequency $\log p$.} We found in realizing this project a surprising relation between mathematical properties and the ideas developed by the composer Olivier Messiaen, as explained in Section \ref{mathmotiv}.

\section{Mathematical Motivation}\label{mathmotiv}
 One of the great unsolved problem in mathematics is the "Riemann Hypothesis" which concerns a formula found by the mathematician Bernhard Riemann relating the prime numbers with the zeros of a mysterious analytic function called the zeta function. My own encounter with this problem came when, while working on quantum statistical mechanics with J. B. Bost, we found a simple quantum mechanical system whose partition function is the Riemann zeta function. For this work I was invited in 1996 to a conference in Seattle celebrating the Norwegian mathematician Atle Selberg's great contributions to the study of the Riemann zeta function. With the help of serendipity I found, when returning from Seattle, that a spectral realization of the zeros of zeta emerged naturally from the quantum physics work I had done with J. B. Bost namely the "BC-system". In the traditional way to look for an interpretation of the zeros of the Riemann zeta function, one looks for them as energy levels of a quantum system (see for instance \cite{BKe}) and what was strange is that in my
   own approach however the zeros appear naturally as "time onsets" \ie in a dual manner (see \cite{Co-zeta}). While this might look perplexing at first sight, I shall briefly explain why this dual point of view is in fact more natural in view of the generalization of the Riemann zeta function in the framework of global fields. Indeed  {  when mathematicians are confronted to a very difficult problem they extend its realm   in order to cast the original unsolved question in a more general context whose larger scope provides ways of considering simpler instances of the problem which, being more tractable, can provide hints of the solution. This is how one meets  much simpler avatars of the Riemann zeta function}. They are associated to a curve $C$ over a finite field $\F_q$. It turns out that these analogues of the Riemann zeta function $\zeta(s)$ are in fact functions of the form $L(q^{-s})$ where $q$ is the cardinality of the finite field over which the curve is defined. Moreover $L(z)$ is a rational fraction and the zeros are determined by the polynomial numerator $P(z)$ whose degree is twice the genus $g$ of the curve. By a famous theorem of Andr\' e Weil all the zeros $z_j$ of $P(z)$ are on the circle of radius $q^{-1/2}$. Thus the zeros of $L(q^{-s})$ are all of real part $\Re(s_k)=\frac 12$ and their imaginary parts are of the form
$$
\Im(s)=\frac{1}{\log q}\left(\alpha_j+2\pi k \right), \ j\in \{1, \ldots 2g\}, \ k\in \Z
$$
where the $-\alpha_j$ are the arguments of the $z_j$ and are determined only modulo $2\pi$. This type of periodic distribution of numbers is very natural as a distribution in ``time". In fact if we specify the choice of the arguments by the fundamental domain $-\pi\leq \alpha_j\leq \pi$ and interpret them as giving {``time onsets"} when notes are played, several striking facts come to the fore:
\begin{enumerate}
\item The functional equation for the $L$ function means that the obtained rhythm is palindromic, \ie with $-\pi\leq \alpha_1\leq \alpha_2\leq \ldots \leq \alpha_g\leq 0\leq \alpha_{g+1}\leq \ldots \leq \alpha_{2g}\leq \pi$ one has
$$
\alpha_{2g+1-j}=-\alpha_j
$$
\item The obtained values are in general irrational numbers.
\item 	The period of the rhythm is given by the  $\frac{2\pi}{\log q}$ and the ``tempo" accelerates when $q$ increases.
\end{enumerate}
The music composer Olivier Messiaen introduced in his work (see \cite{Messiaen0} and the full treatise \cite{Messiaen}) both the palindromic feature (which he called "rythmes non-r\' etrogradables") and the irrational {"time onsets"} of attack. It is worth drawing a picture of the type of palindromic rhythm which is delivered by the zeros of the $L$-function of a curve over a finite field, and we shall undertake doing that below.

\section{Explicit examples of motivic rhythms}
 In this section we exhibit in visual form the motivic rhythms of genus $g=5$ in order to illustrate the general construction. The choice of the number $g=5$ could be replaced by any other number, the number of time onsets is $2g=10$. The motives are determined by choosing a curve $C$ of genus $g=5$ defined over the field of rational numbers and then reducing the curve modulo each of the prime numbers $2,3,5,7,11, \ldots$ except for those finitely many primes called "of bad reduction" for which the reduction is ill behaved. The computationally most convenient curves are the "hyperelliptic" ones. They are called in this way because the simplest of them (those with $g=1$) are the well known elliptic curves which are a jewel of algebraic geometry. An hyperelliptic curve over $\Q$ is given by an equation of the form $y^2=P(x)$ where $P(x)$ is a polynomial with integer coefficients and whose degree determines the genus so that we shall deal with degree $11$ to get a curve of genus $5$. The primes of bad reduction are those which  divide the discriminant of $P(x)$ and thus we shall pick the simplest ones whose discriminant does not involve any prime between $7$ and $67$ so that we can then explore the motives associated to this sequence of primes.

Starting with such an hyperelliptic curve  one obtains for each prime not dividing the discriminant of $P$ a curve over the finite field $\F_p$ and one can compute the zeros of its $L$-function. Thus the choice of the polynomial $P(x)$ is dictated by the absence of a prime divisor of the discriminant between $7$ and $67$. We start with
$$
C_1: \ y^2=x^{11}-4 x^{10}+15 x^8-40 x^6+20 x^5+25 x^3-25
$$
We give below the list of polynomials of degree $10$, associated to the primes $7,11,\ldots, 67$ and whose zeros give the eigenvalues of the Frobenius acting on $H^1(C_1/p)$ where $C_1/p$ is the reduction of the hyperelliptic curve $C_1$ modulo the prime $p$. The degree of these polynomials is twice the genus, \ie $5\times 2=10$ and the term of degree $10$ is always of the form $p^5x^{10}$. The complex zeros are all of modulus $p^{-1/2}$ and the arguments of the zeros give us palindromic rhythms which are fixing the interpretation associated to $C_1/p$.

\begin{figure}[H]
\begin{center}
\includegraphics[scale=0.6]{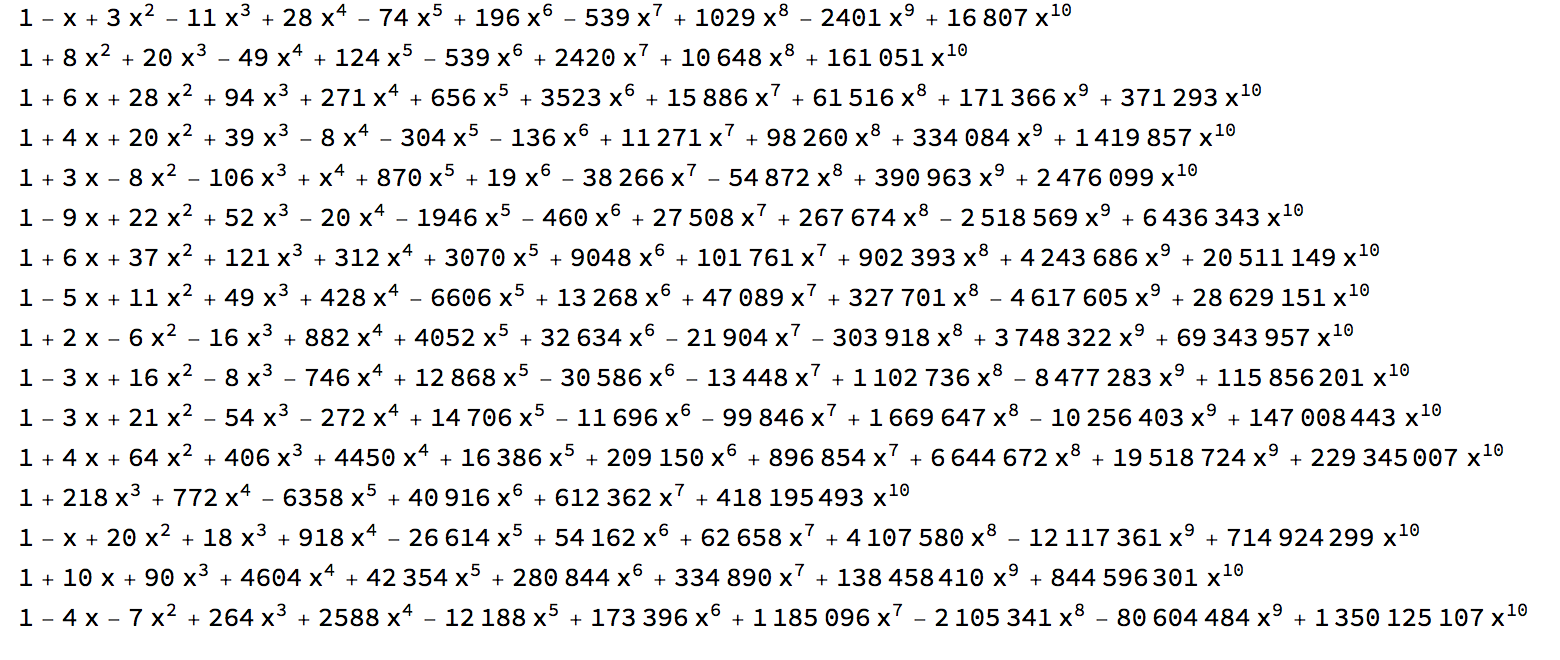}
\end{center}
\caption{Polynomials associated to $H^1(C_1/p)$ for $7\leq p\leq 67$. \label{poly1} }
\end{figure}
{ We did not display the prime $p$, $7\leq p \leq 67$, corresponding to each of these polynomials but one sees that the coefficients of $x^{10}$ increase when going down the list as they should since they are given by the fifth power $p^5$. The "palendromic" property of these polynomials is, with $P(z)$ associated to the prime $p$ the equality
$$
p^5z^{10}P(1/(pz))=P(z)
$$
which corresponds to the functional equation of the zeta function \ie the replacement $s\mapsto 1-s$ whose effect on $z=p^{-s}$ is
$$
z=p^{-s}\mapsto p^{-(1-s)}=1/(pz).
$$
}
Since the rhythms repeat periodically we just show them for one period and for simplicity we do not show the acceleration of the tempo. { We first display them linearly and then in a circular manner which fits with their periodic property and is closest to their mathematical meaning.}
\begin{figure}[H]
\begin{center}
\includegraphics[scale=0.75]{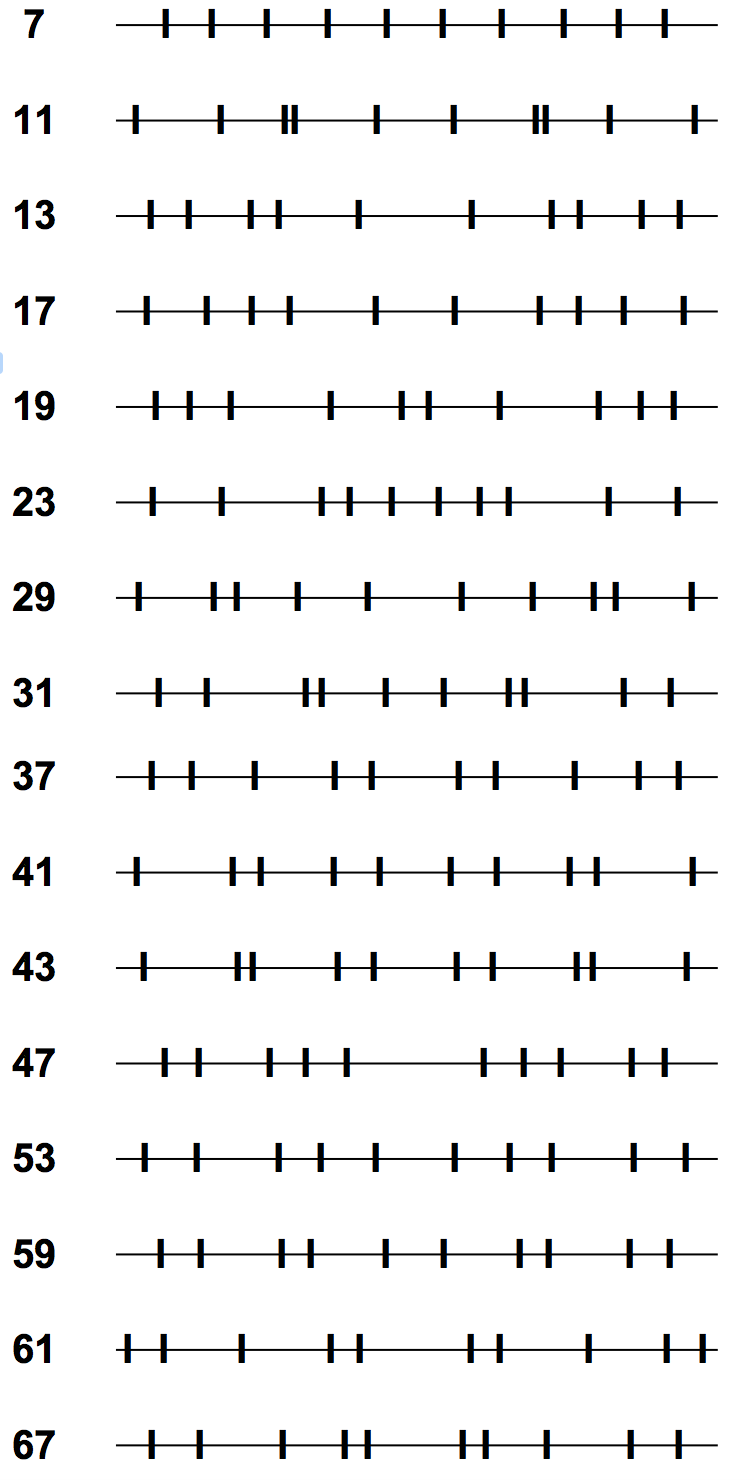}
\end{center}
\caption{Motivic Rhythms for $H^1(C_1/p)$, up to $p=67$.  \label{rhythm1} }
\end{figure}

\begin{figure}[H]
\begin{center}
\includegraphics[scale=0.3]{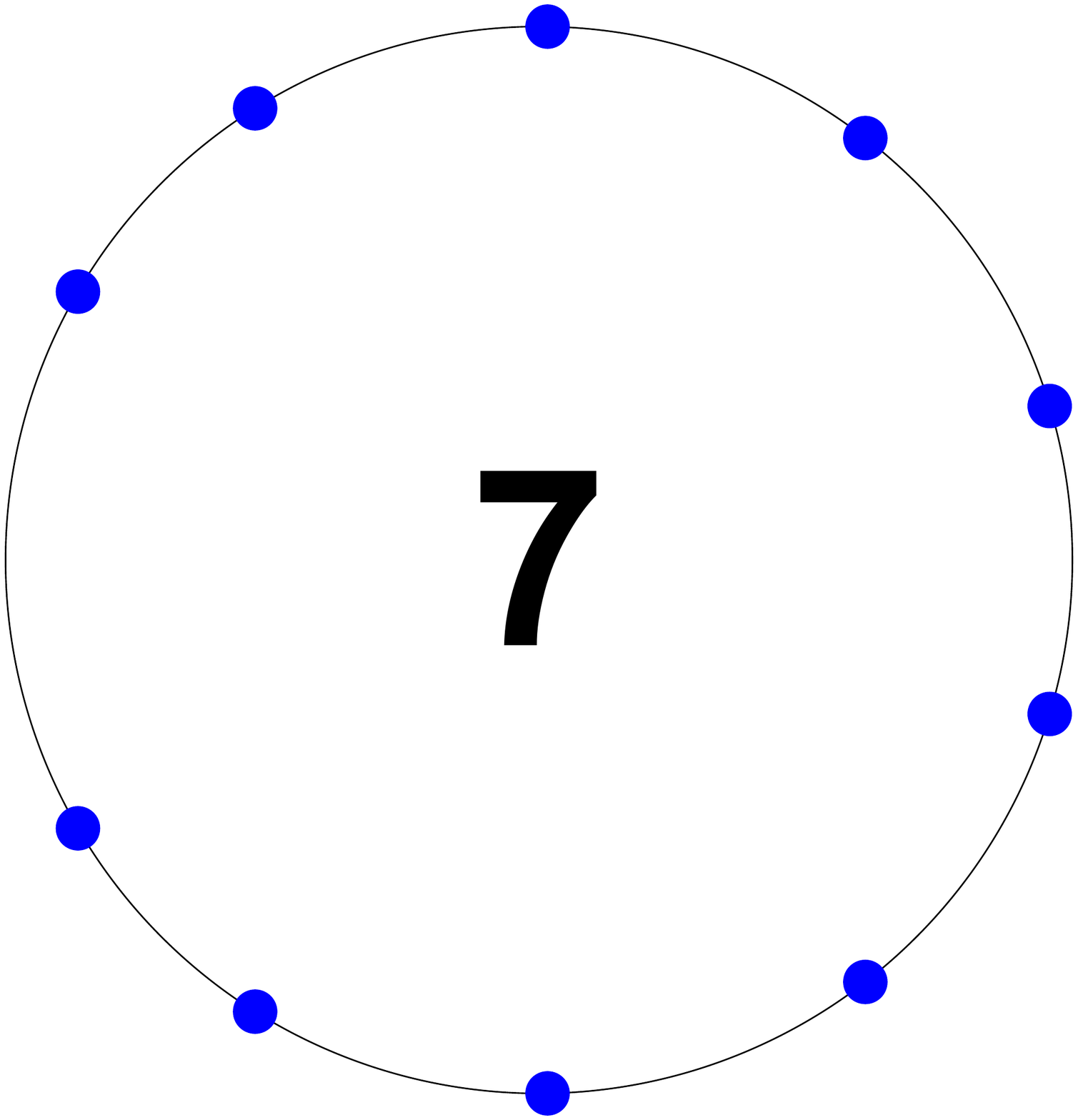}
\end{center}
\caption{Motivic Rhythm for $H^1(C_1/p)$, $p=7$.\label{v1} }
\end{figure}
\begin{figure}[H]
\begin{center}
\includegraphics[scale=0.3]{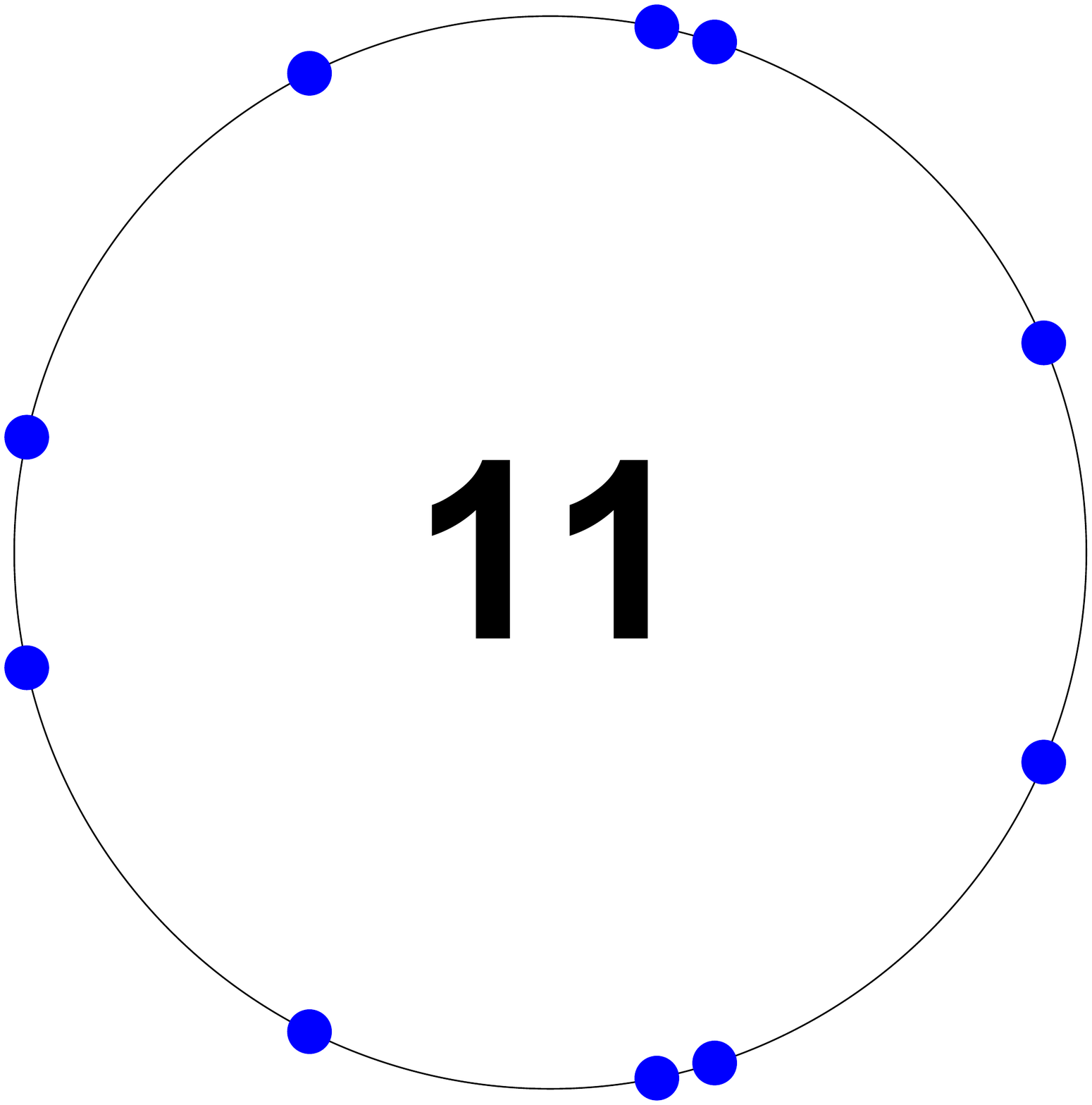}
\end{center}
\caption{Motivic Rhythm  for $H^1(C_1/p)$, $p=11$. \label{v2} }
\end{figure}
\begin{figure}[H]
\begin{center}
\includegraphics[scale=0.3]{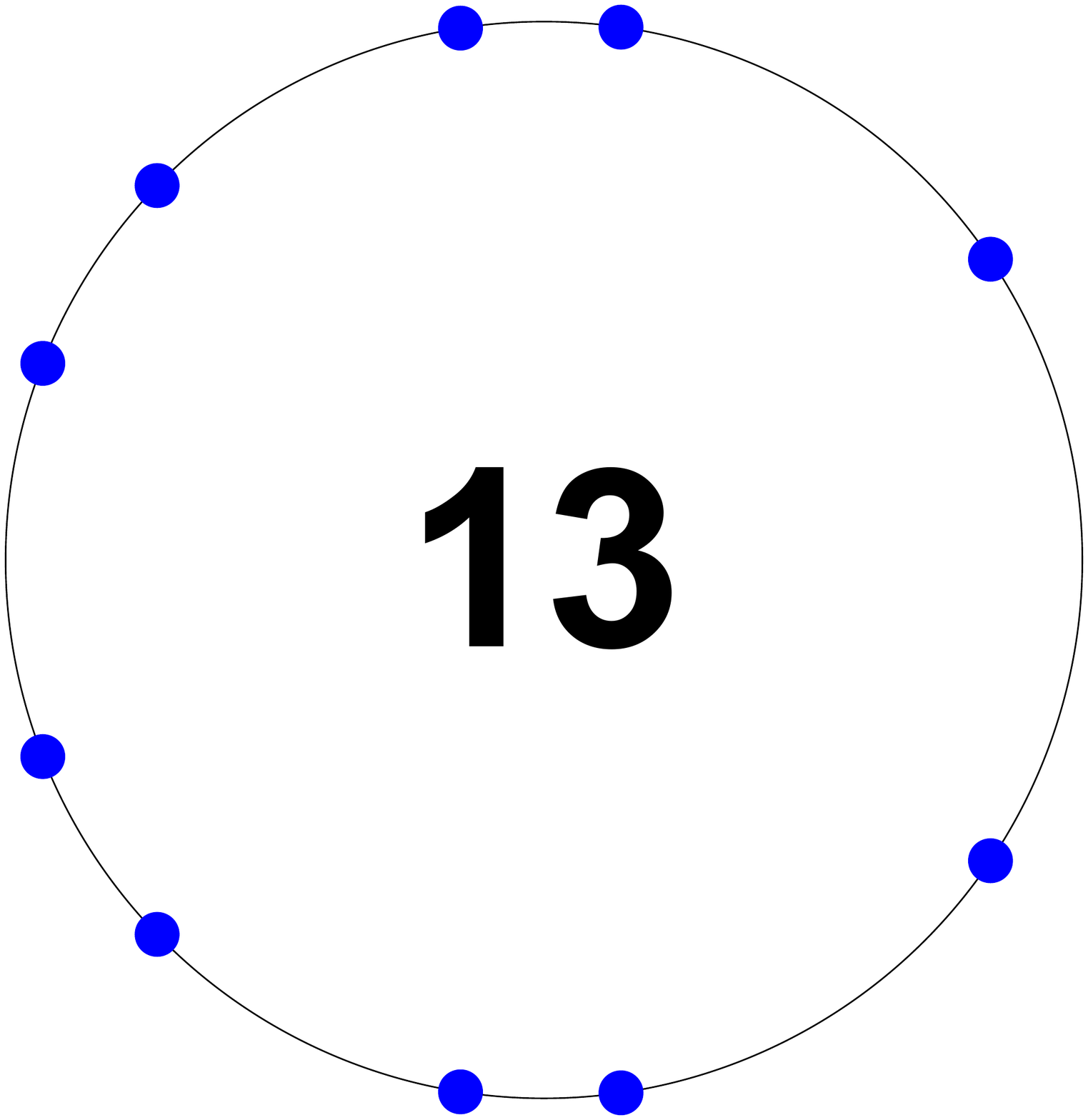}
\end{center}
\caption{Motivic Rhythm  for $H^1(C_1/p)$, $p=13$.\label{v3} }
\end{figure}
\begin{figure}[H]
\begin{center}
\includegraphics[scale=0.3]{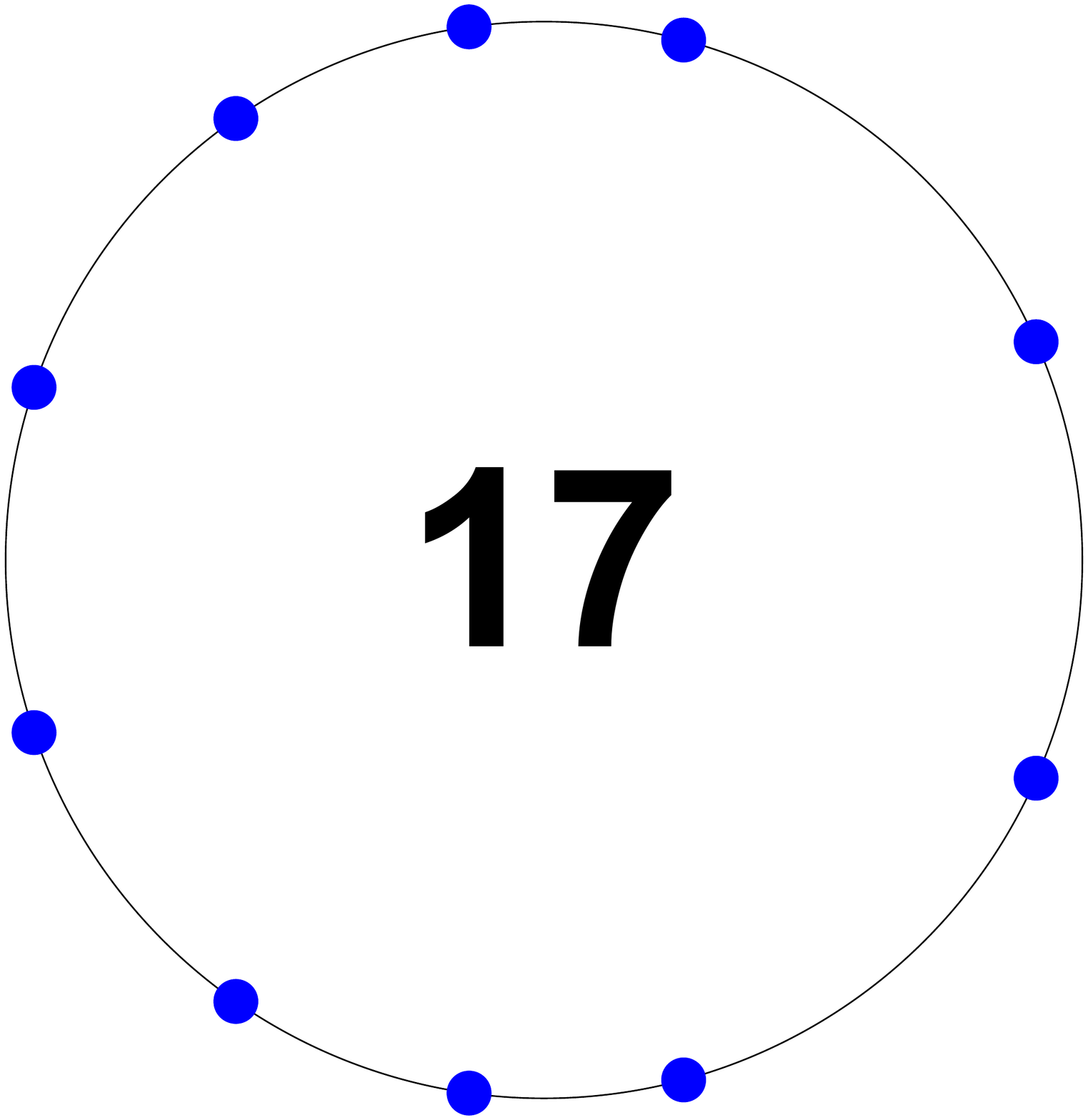}
\end{center}
\caption{Motivic Rhythm  for $H^1(C_1/p)$, $p=17$. \label{v4} }
\end{figure}
\begin{figure}[H]
\begin{center}
\includegraphics[scale=0.3]{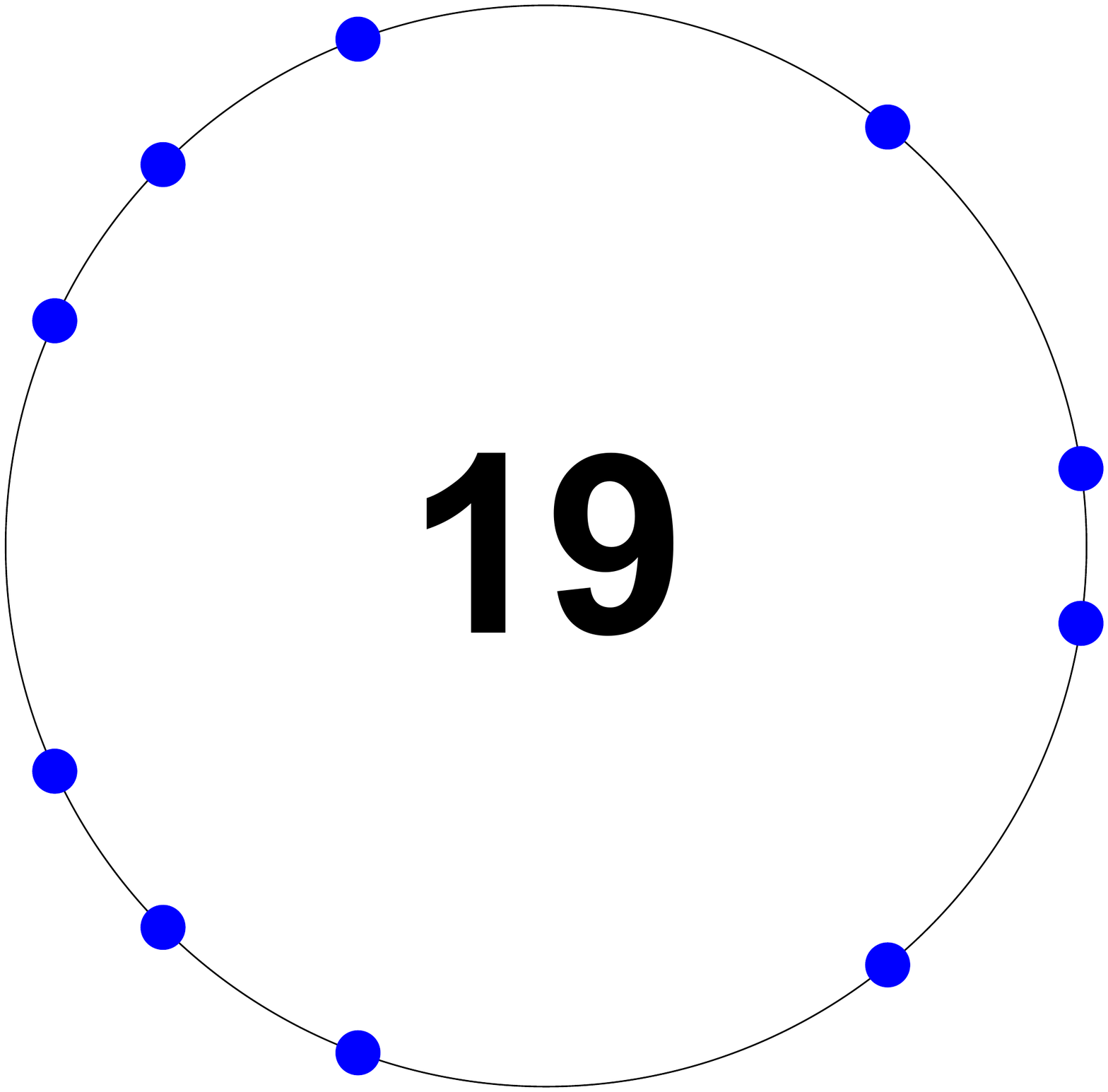}
\end{center}
\caption{Motivic Rhythm  for $H^1(C_1/p)$, $p=19$. \label{v5} }
\end{figure}

We give the next example of hyperelliptic curve of genus $5$ with good reduction at all primes between $7$ and $67$. It is given by the equation:

 $$
C_2: \  y^2= x^{11}-60 x^7-64 x^6-320 x^4-380 x^3-512 x-640
 $$
 \begin{figure}[H]
\begin{center}
\includegraphics[scale=0.6]{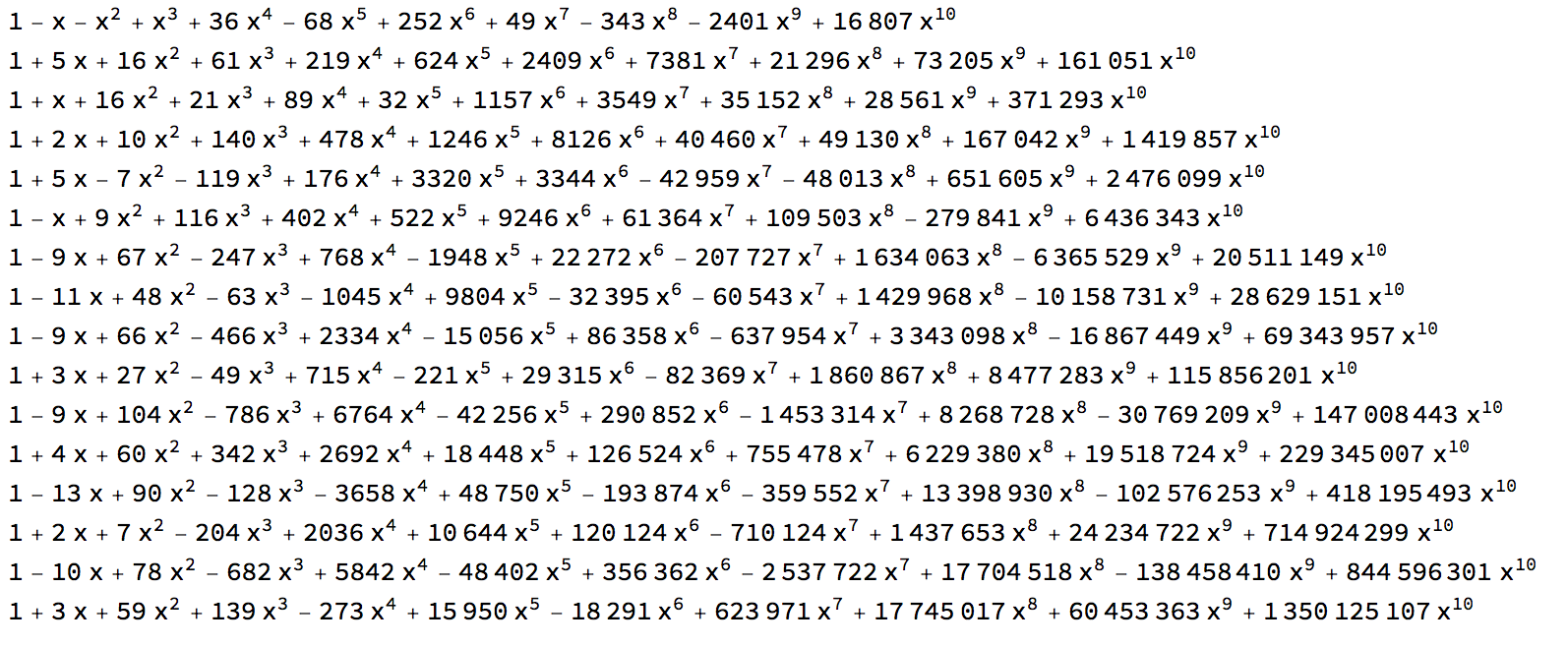}
\end{center}
\caption{Polynomials associated to $H^1(C_2/p)$ for $7\leq p\leq 67$. \label{poly1} }
\end{figure}

 \begin{figure}[H]
\begin{center}
\includegraphics[scale=0.75]{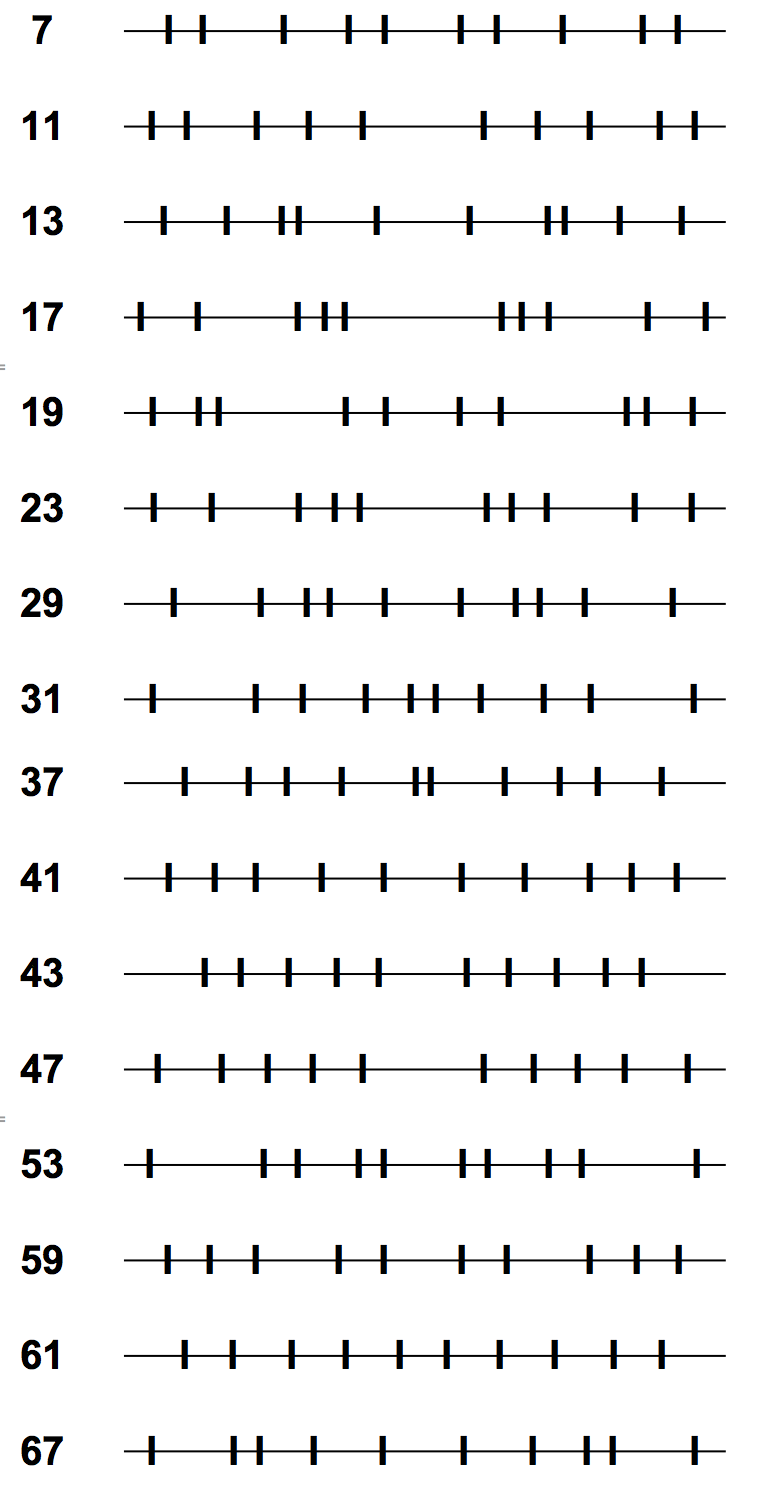}
\end{center}
\caption{Rhythms associated to $H^1(C_2/p)$ for $7\leq p\leq 67$. \label{poly1} }
\end{figure}

\begin{figure}[H]
\begin{center}
\includegraphics[scale=0.25]{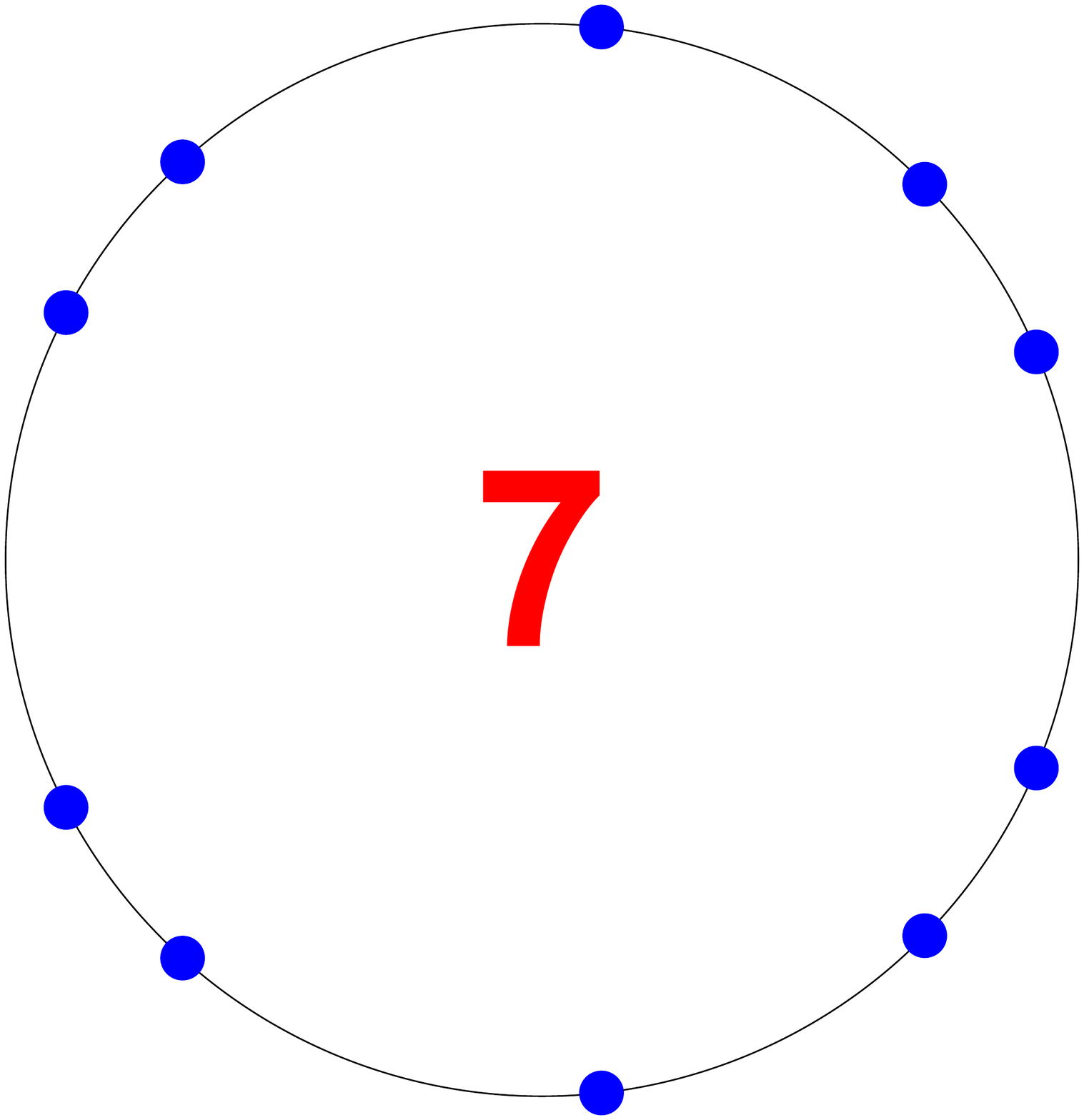}
\end{center}
\caption{Motivic Rhythm for $H^1(C_2/p)$, $p=7$ \label{w1} }
\end{figure}
\begin{figure}[H]
\begin{center}
\includegraphics[scale=0.23]{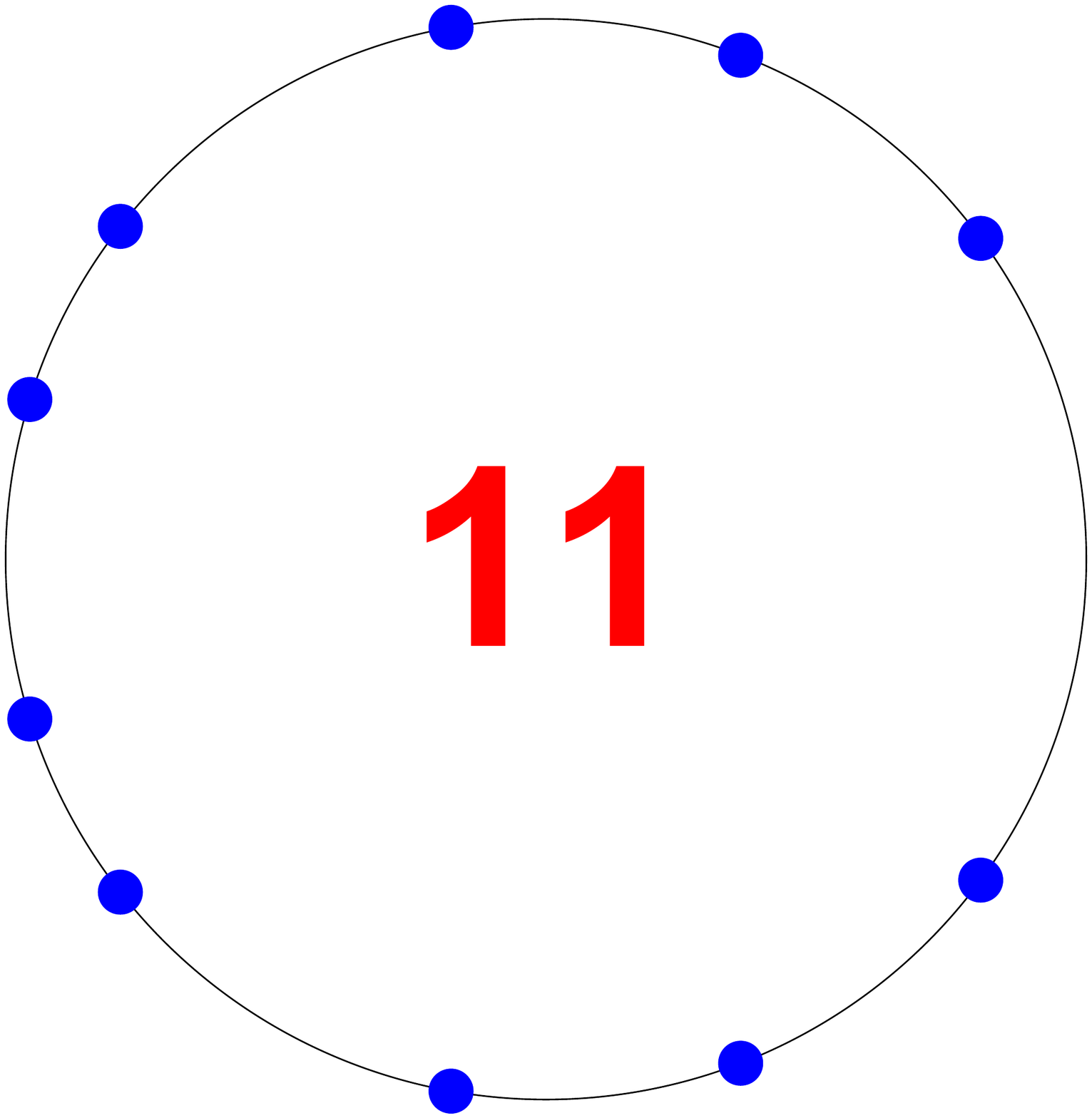}
\end{center}
\caption{Motivic Rhythm  for $H^1(C_2/p)$, $p=11$\label{w2} }
\end{figure}
\begin{figure}[H]
\begin{center}
\includegraphics[scale=0.20]{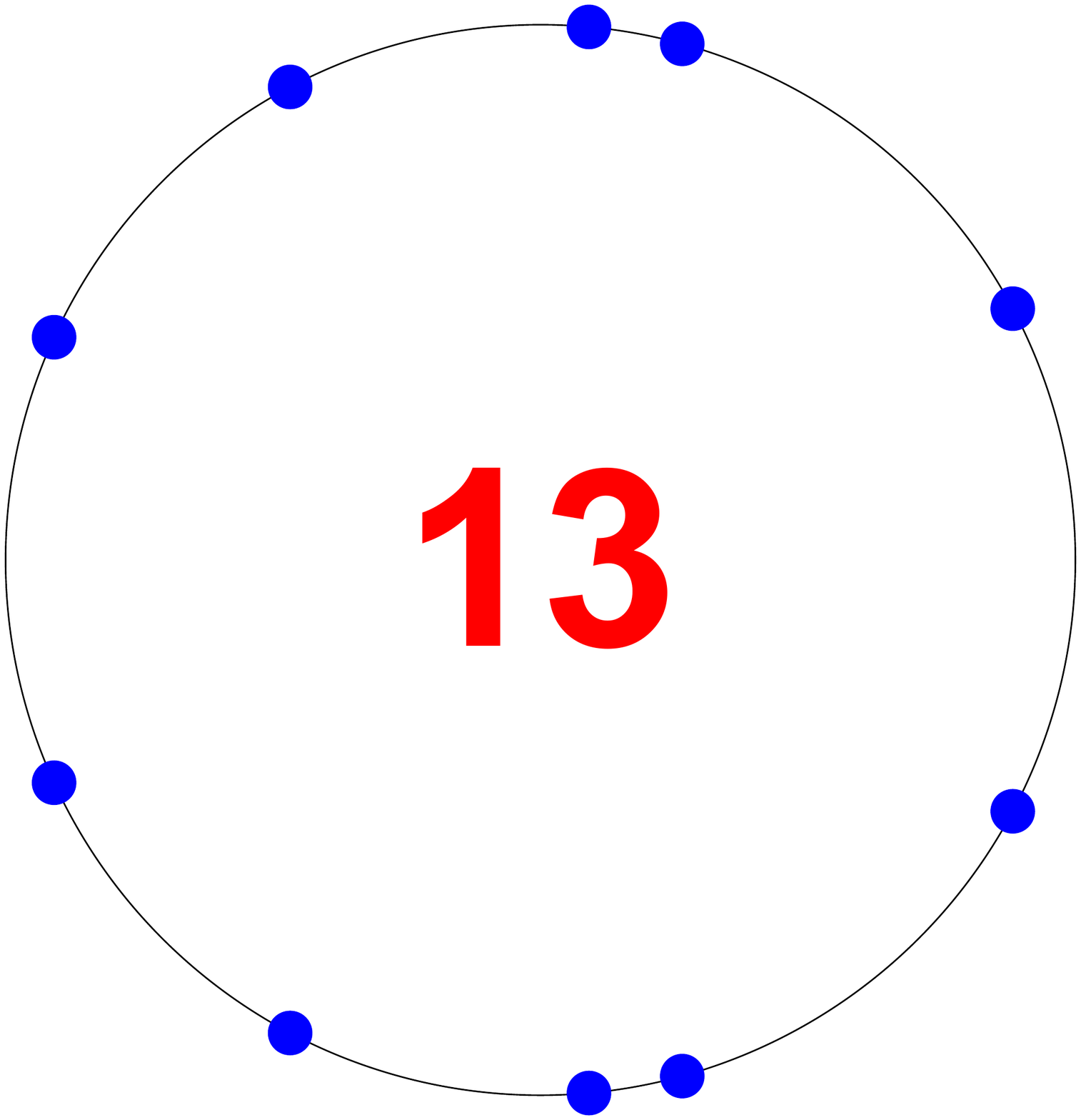}
\end{center}
\caption{Motivic Rhythm  for $H^1(C_2/p)$, $p=13$ \label{w3} }
\end{figure}
\begin{figure}[H]
\begin{center}
\includegraphics[scale=0.17]{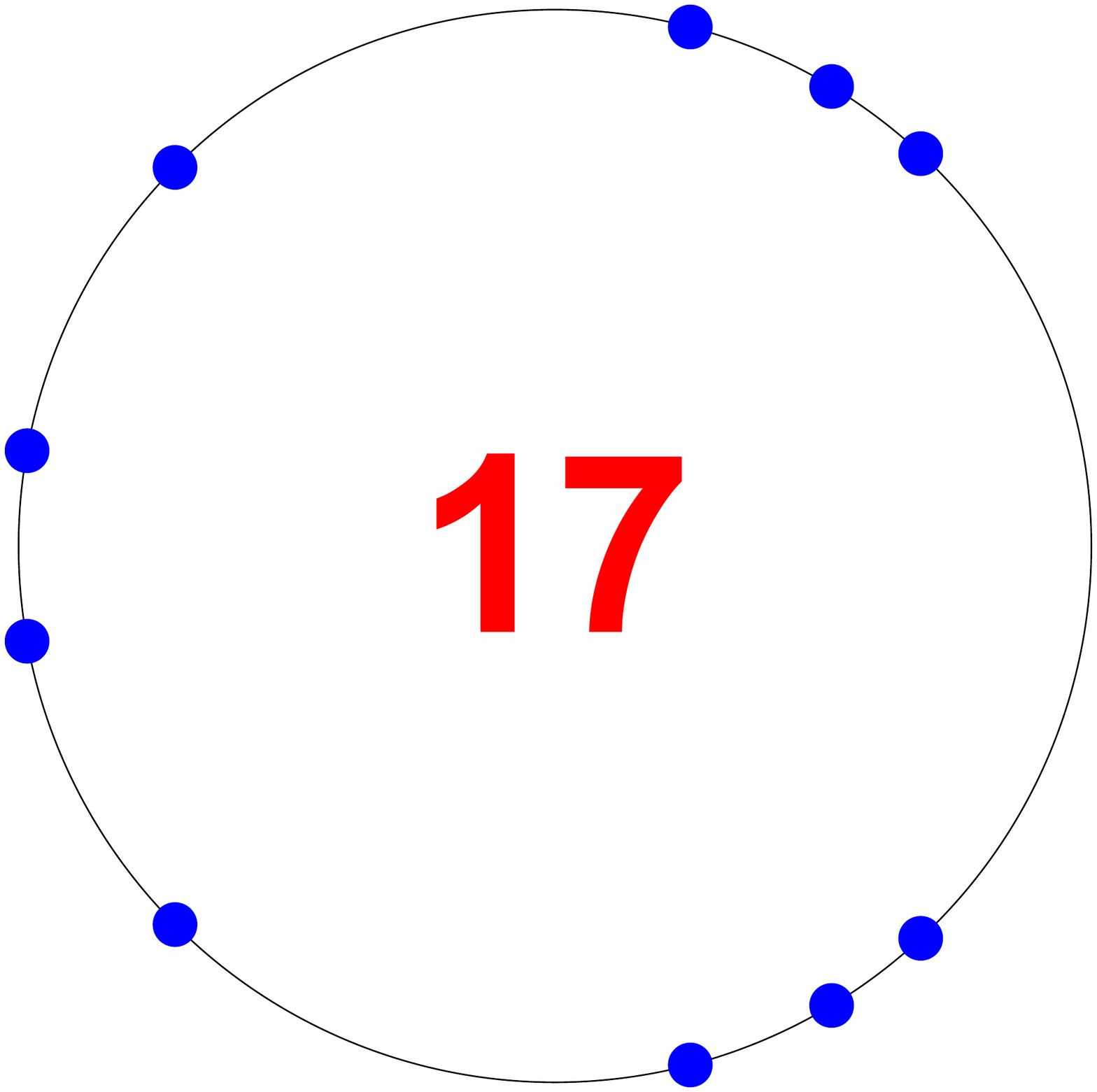}
\end{center}
\caption{Motivic Rhythm  for $H^1(C_2/p)$, $p=17$ \label{w4} }
\end{figure}
\begin{figure}[H]
\begin{center}
\includegraphics[scale=0.15]{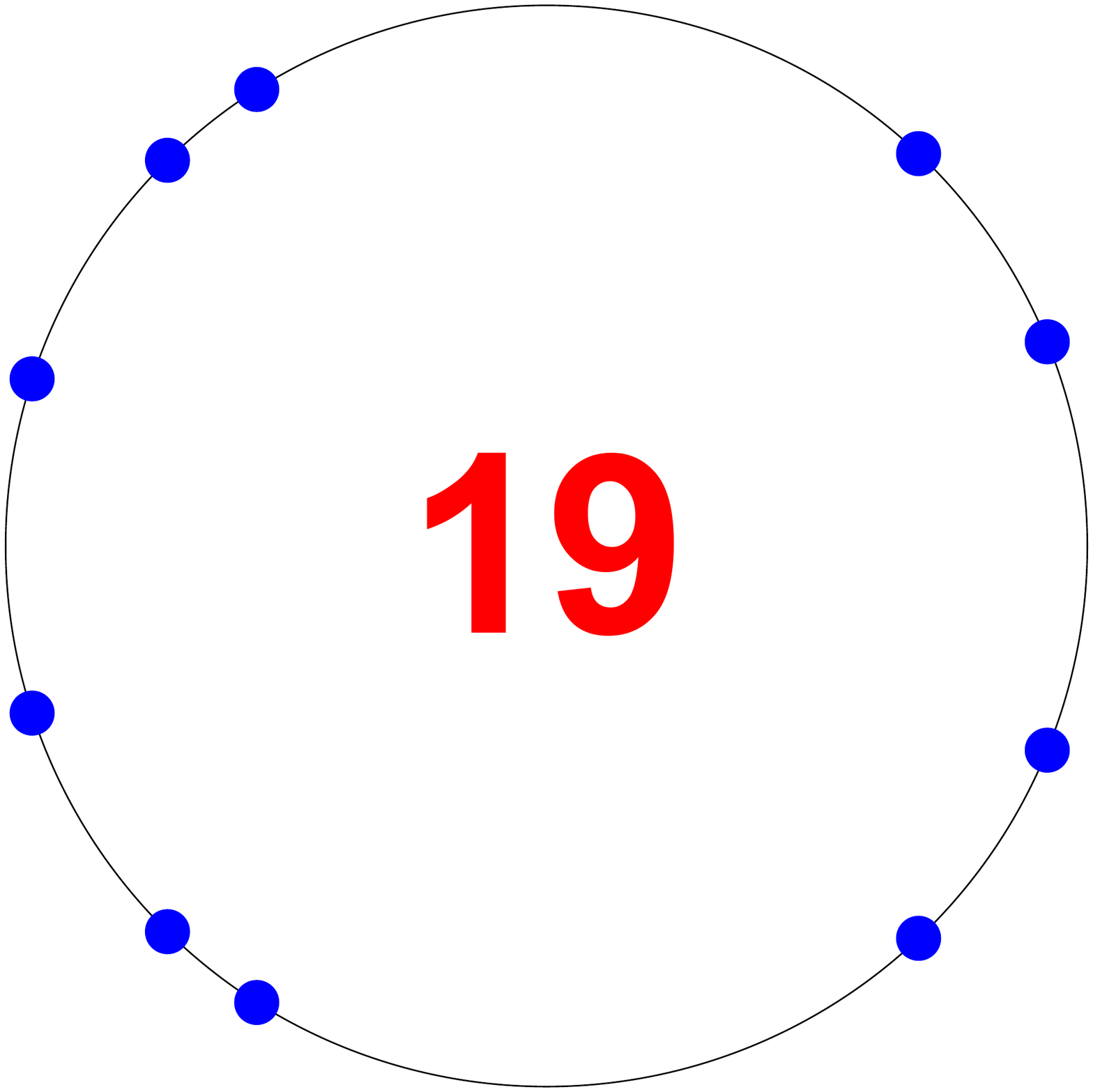}
\end{center}
\caption{Motivic Rhythm  for $H^1(C_2/p)$, $p=19$.\label{w5} }
\end{figure}
   In the next section we explain how we devised a mathematical formula for a basic score attached to the primes between $7$ and $67$.

    There are remarkable interrelations, both cognitive and theoretical, between pitches and rhythms which are well documented in the literature (see \cite{Pressing}, \cite{Rhan}). We have chosen to keep the rhythmic part totally independent of the  choice of a collection of pitches associated to primes discussed in the next section. The first  reason to maintain this independence is that the motivic rhythms apply to the playing of any score with the right number of notes and palendromic property. The second reason is to stay in accordance with the mathematical origin of the time onsets which would be lost if they were reinterpreted in terms of frequencies.

 \section{Basic score}\label{musicpiece}

  The motivic rhythms coming from a curve $C$ depend on a prime number and since one has the freedom to choose the pitches of the palindromic song which will be rhythmically interpreted by the motive, one finds that one needs to associate a song \ie here an ordered collection of $5$ notes (which will yield a total of $10$ notes by symmetry)  to each prime number (here those between $7$ and $67$). What matters is that the notes associated to the primes  allow one to recognize which prime it is. One may try at first simple methods like the name of the prime in a given language but this is obviously far from "canonical". To obtain a canonical process one needs to use the information coming from the well-known  mathematics behind the frequencies in the usual musical scale which are given by the integral powers of the number $q:=2^{\frac{1}{12}}$ (which is numerically close to $3^{\frac{1}{19}}$).
 In this section we describe how we obtain from a mathematical formula a palindromic score of ten notes for each prime between $7$ and $67$. The idea is that the score associated to the prime $p$ should be related to the number $\sqrt p$ because the motive is of weight $\frac 12$. { This weight of the motive corresponds to the shrinking size of the circle where the  zeros are as shown in Figures \ref{w1}... \ref{w5}. The inverse radius of this circle is $\sqrt p$.  Moreover since the musical scale is the set $\{q^n\}$ of integral powers of the number $q:=2^{\frac{1}{12}}$ it is natural to start by the integral power $q^{n_1(p)}$ which is the closest to $\sqrt p$. This corresponds to a pitch naturally associated to $p$. Thus as  starter we let}
 $$
 n_1(p):= {\rm IntegerPart}(6\log p/\log 2)
 $$
 Since { in order to build the next notes in the score } we need to obtain integers and to make sure that they suffice to characterize $p$, we have chosen to use the next terms $n_j(p)$ in the  continued fraction expansion  $$6\log p/\log 2=n_1+1/(n_2+1/(n_3+1/(n_4+1/(n_5+\ldots))))$$ To respect the respective roles of the $n_j(p)$ it is natural to convert them to the following list
 $$
 (n_1(p),n_1(p)-n_2(p),n_1(p)-n_2(p)+n_3(p),n_1(p)-n_2(p)+n_3(p)-n_4(p),$$ $$n_1(p)-n_2(p)+n_3(p)-n_4(p)+n_5(p))
 $$
 This provides us with a score which we complete by symmetry, to make it palindromic,  in order to be able to reflect the palindromic rhythms,  for each prime between $7$ and $67$.

 \begin{figure}[H]
\begin{center}
\includegraphics[scale=0.95]{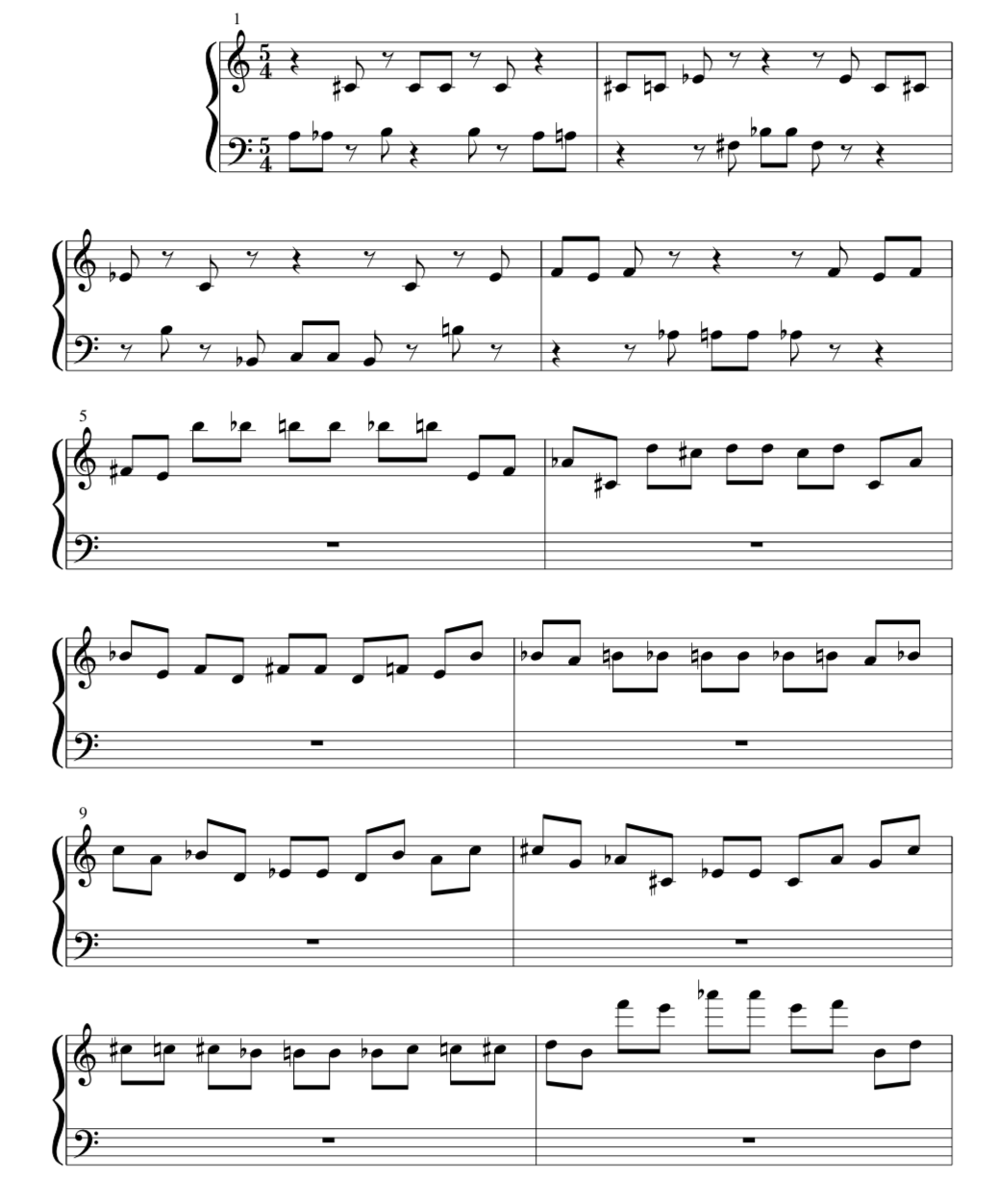}
\end{center}
\caption{Simple musical line associated to $7\leq p\leq 47$.\label{start} }
\end{figure}

 \begin{figure}[H]
\begin{center}
\includegraphics[scale=0.95]{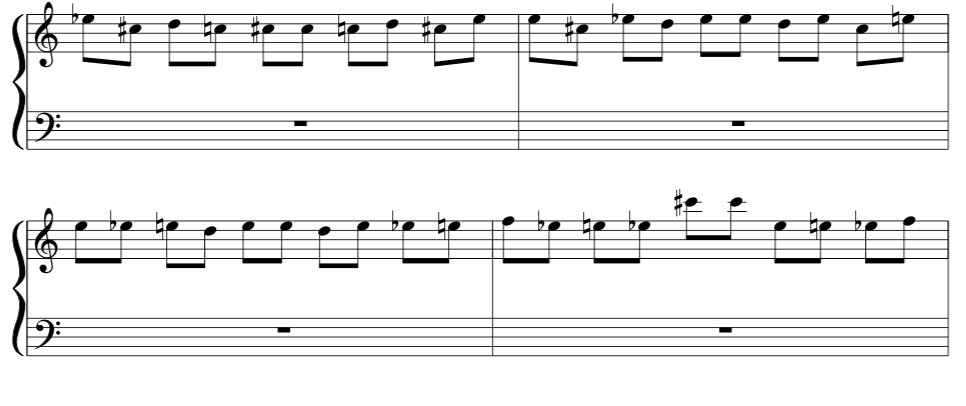}
\end{center}
\caption{Simple musical line associated to $53\leq p\leq 67$.\label{start} }
\end{figure}

 There remains at this point the musical task to harmonize the above simple score and there is a clear freedom in doing that. One constraint is that the chords ought to be played at the same time as the notes of the initial tune in order to let them undergo the change of rhythm dictated by the ``motivic interpretation". {  We have devised  a simple one used in the creation of the video, but we put it as a challenge to musicians to devise a really good one!}
{ Without this constraint one can produce an harmonization which is more dynamical and that allows one to get familiar with the notes of the initial tune (the tunes of the first two primes $7$ and $11$ are repeated twice) it can be downloaded  at \href{https://www.dropbox.com/s/df3vvuhkvel5eap/musescoremi.mp3?dl=0}{\color{blue}link to download the mp3}. We shall use it also below in the visual illustration.}

 \section{Visual illustration:  { The dance of primes} }

 In order to illustrate visually the succession of primes and then the playing of each motivic rhythm we have created a small choreography.

 \subsection{The first minute: the Eratosthenes sieve}
\begin{figure}[H]
\begin{center}
\includegraphics[scale=0.8]{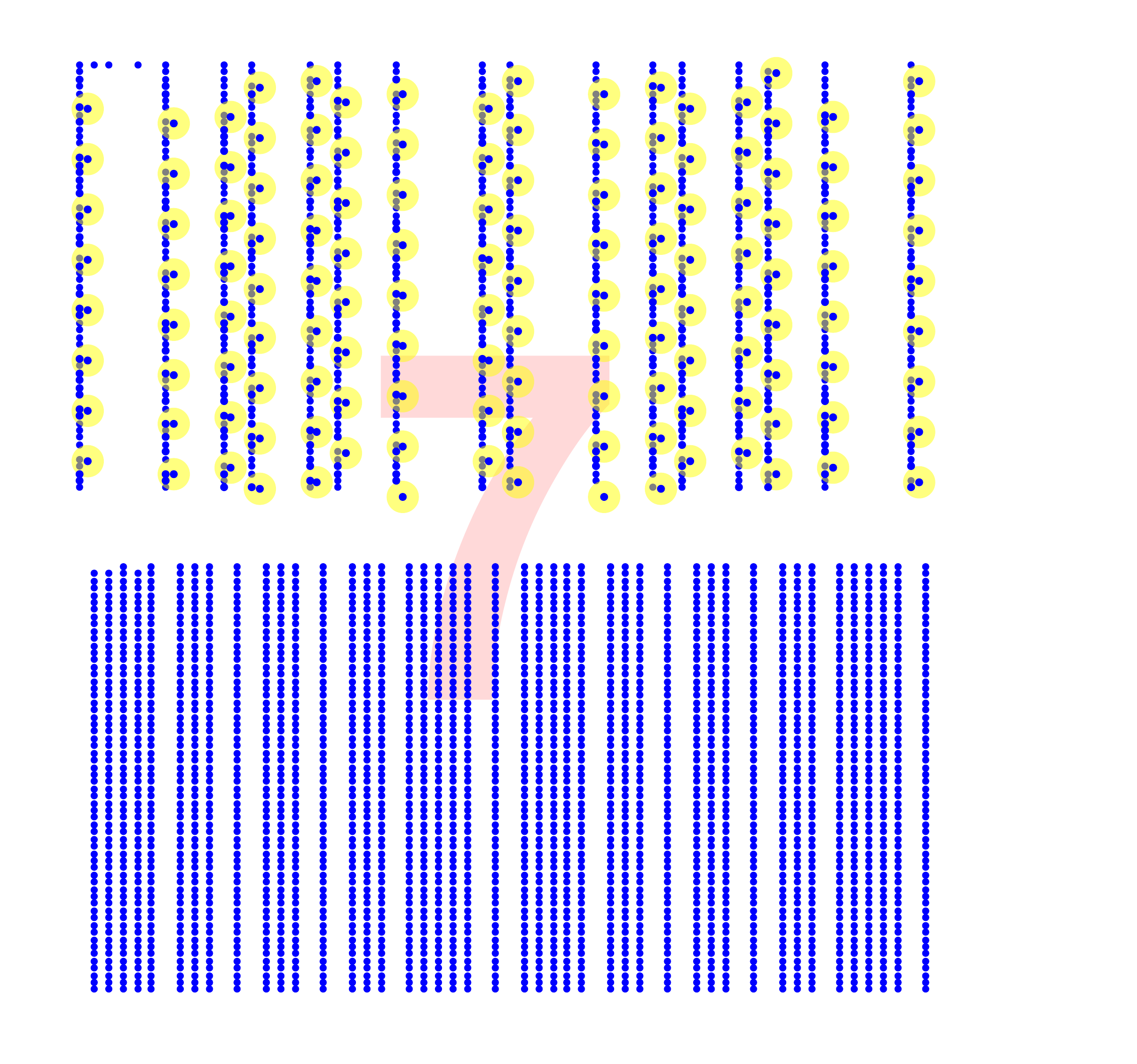}
\end{center}
\caption{All multiples of $7$ except $7$ itself are descending  in the rectangle below. This shows the first step when the dancers corresponding to multiples of $7$ make a side step to the right. \label{7} }
\end{figure}

\begin{figure}[H]
\begin{center}
\includegraphics[scale=0.8]{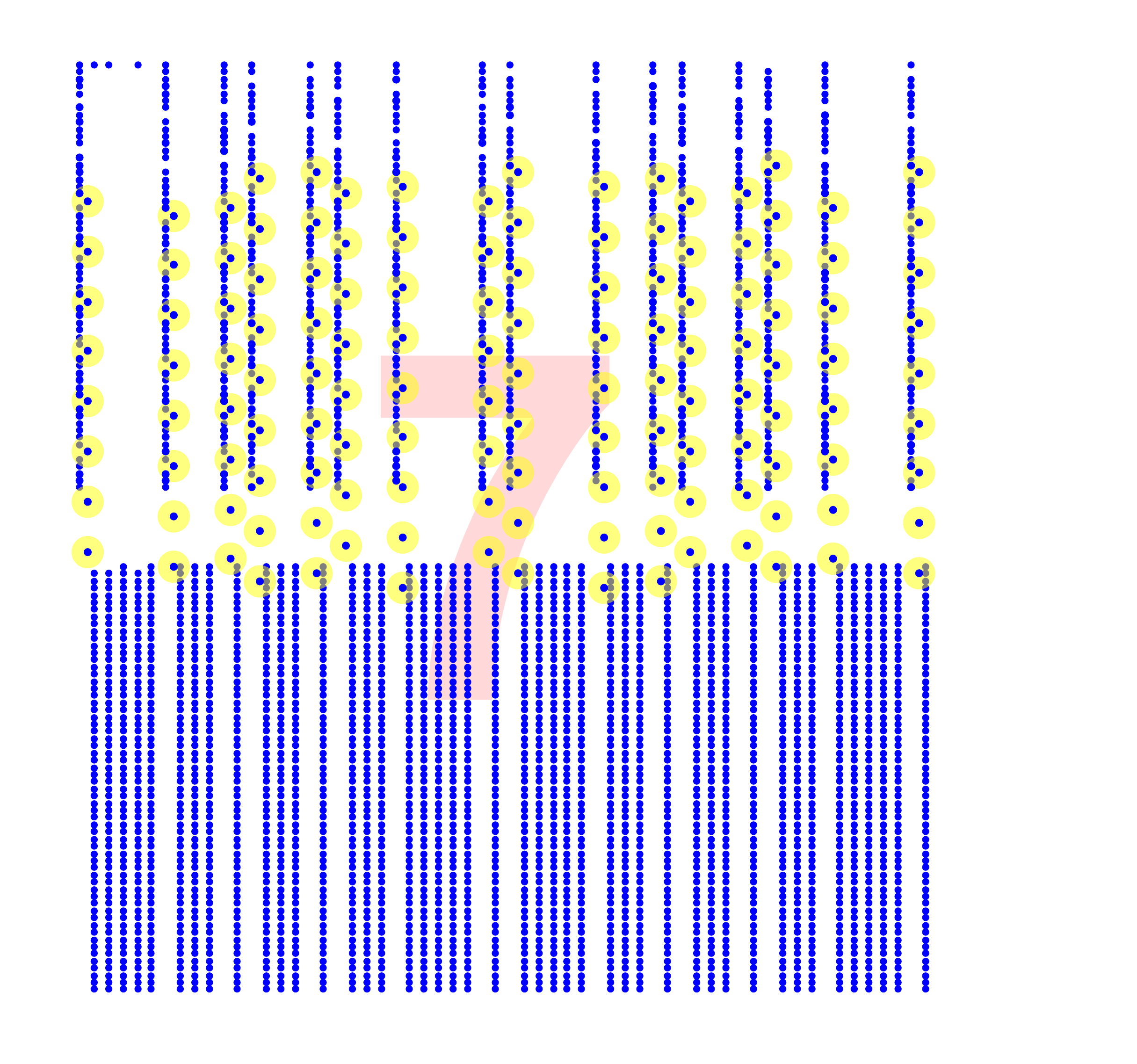}
\end{center}
\caption{All multiples of $7$ except $7$ itself are descending  in the rectangle below. { While this movement takes place the little tune corresponding to $7$ is played}. \label{7half} }
\end{figure}

\begin{figure}[H]
\begin{center}
\includegraphics[scale=0.8]{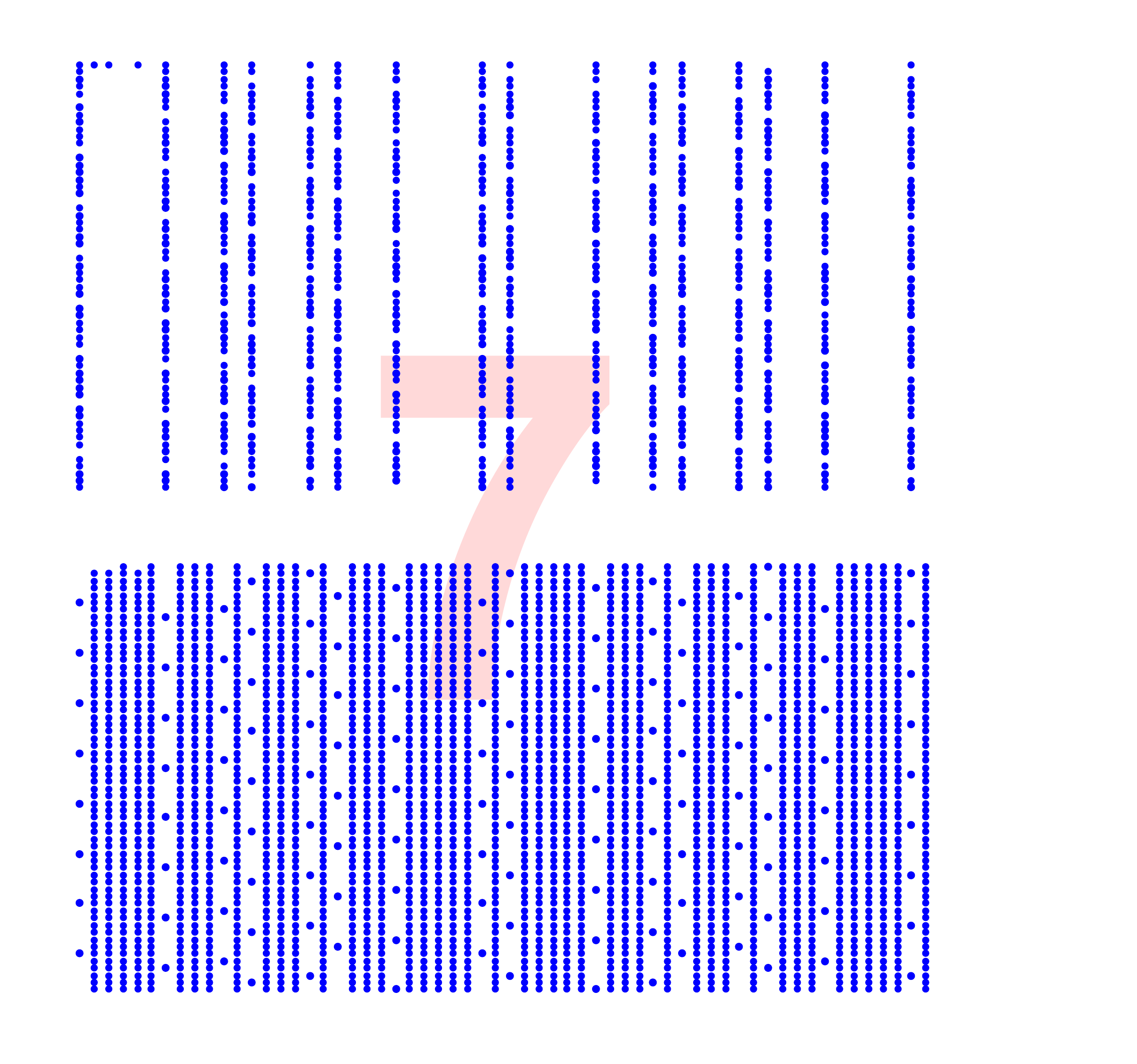}
\end{center}
\caption{All multiples of $7$ except $7$ itself have now descended and found their place  in the rectangle below.{ While this movement took place the music corresponding to $7$ was played}. \label{7bis} }
\end{figure}

  In the first part the music is played with equally spaced time onsets  and the visual illustration is a choreography of the Eratosthene sieve which is devised as follows. There are two rectangles as in Figure  \ref{start}. Each of these rectangles contains $60 \times 60=3600$ possible positions { each of which corresponds to a number between $1$ and $3600$. The first line for instance corresponds to numbers between $1$ (on the top left) and $60$ (on the top right), the second line to numbers between $61$ and $120$ and so on.  At the start all the multiples  of $2$ (\ie all even numbers) except for $2$ itself, have been brought down to the lower rectangle. The same operation has been performed for all multiples of $3$ and all multiples of $5$. Being a multiple of one of these three first primes is unaffected if one adds a multiple of $60$ and thus the obtained geometric picture is formed of a union of vertical lines as shown  in Figure  \ref{start}. The only vertical line in the initial picture which does not correspond to a prime $p<60$, is the vertical line associated to $49=7^2$. All other non-prime have been eliminated as multiples of $2$, $3$ or $5$. The choreography starts when the music associated to $7$ is played, the dancers in the above rectangle which are at a position which is a multiple of seven (except at  seven itself) all descend in the lower rectangle as shown in Figures \ref{7}, \ref{7half}, \ref{7bis} and find their place there. Next when the music of $11$ is played all dancers at positions which are multiples of eleven descend except eleven itself. This continues for all prime numbers until $67$.} At the end of the choreography the numbers which remain in the upper rectangle are all the prime numbers between $2$ and $3600$ (the number $1$ is ignored). The number of dancers required is about 460 because at the start there are only $963$ dancers in the upper rectangle while at the end only the primes less than $3600$ remain and there are $503$ such primes. The video is available by clicking on the following link:
 \href{https://www.dropbox.com/s/ld1z64em5u3mxjx/fullvideo.mp4?dl=0}{\color{blue}Video}. { For the "dynamical" harmonization of the equally played ``dance of primes" one can click on the following link:  \href{https://www.dropbox.com/s/fvf7rctfpl3eycw/animation-equal-good.mp4?dl=0}{\color{blue}Dance of Primes}. Note that the song of each prime is here repeated twice with dancers standing still during the repetition.}

 \begin{figure}[H]
\begin{center}
\includegraphics[scale=1]{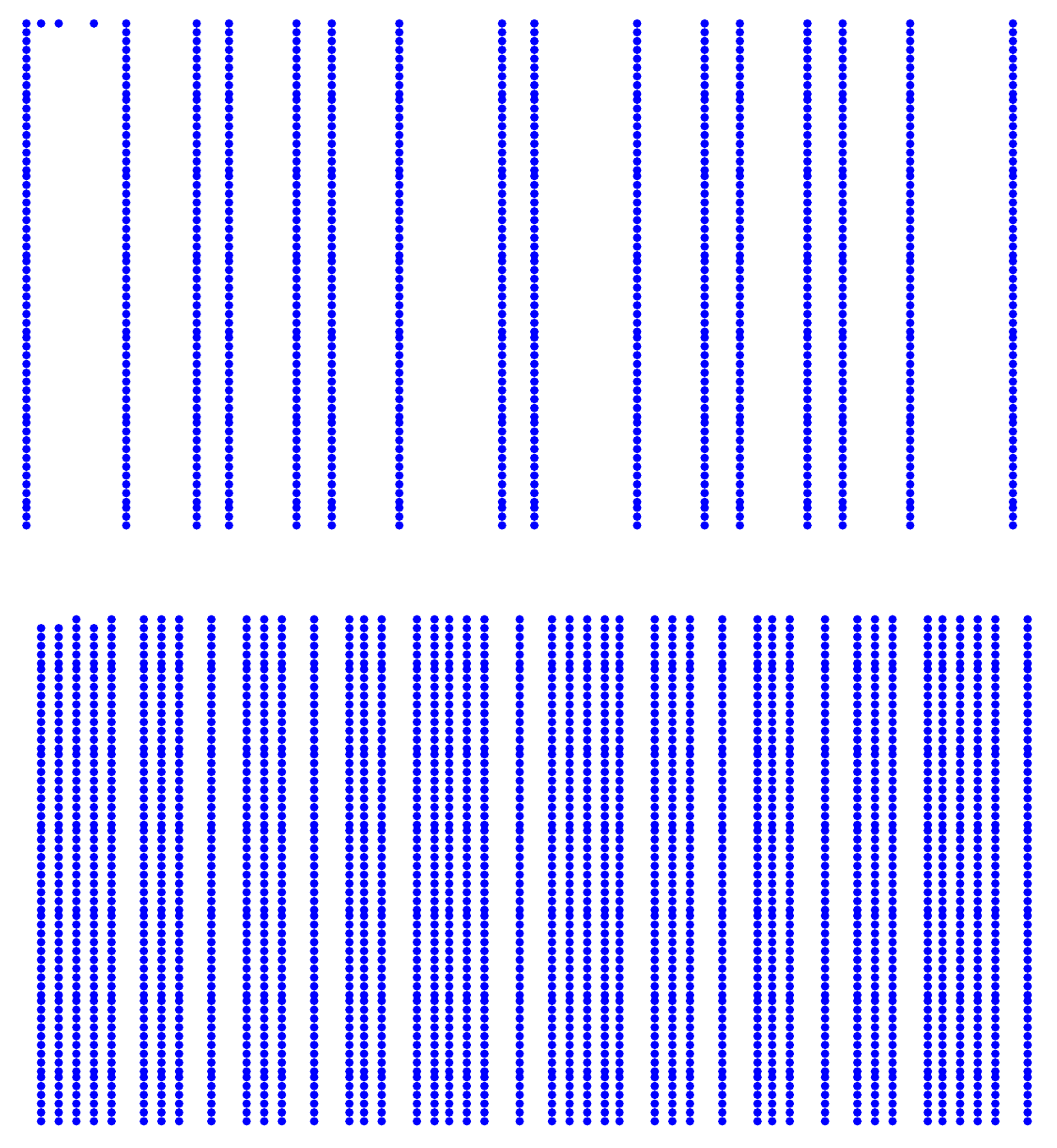}
\end{center}
\caption{Starting configuration: all multiples of $2$ except $2$ itself, and similarly for $3$ and $5$ have been brought down in the rectangle below.\label{start} }
\end{figure}
\begin{figure}[H]
\begin{center}
\includegraphics[scale=0.8]{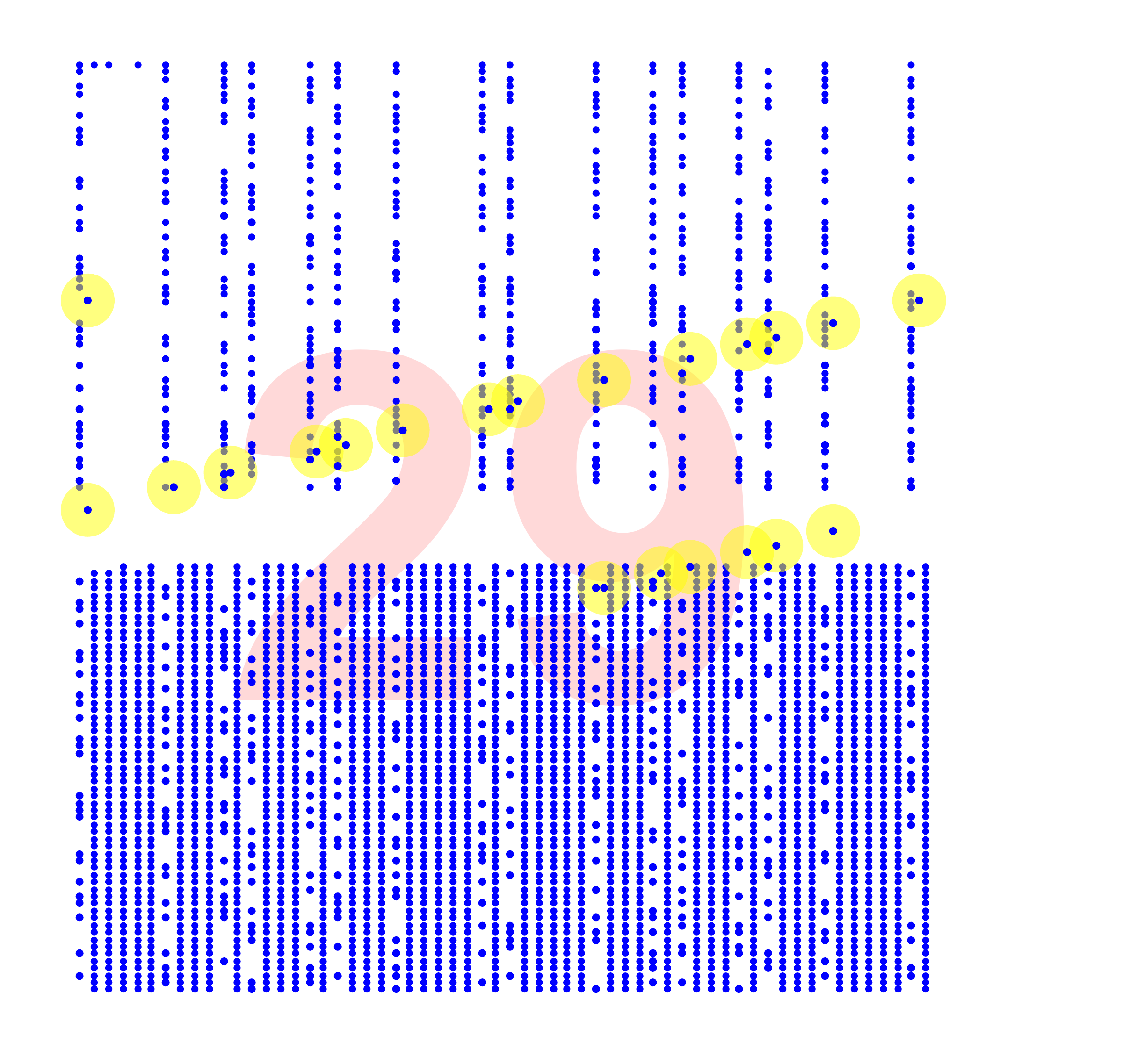}
\end{center}
\caption{All multiples of $29$ except $29$ itself are descending  in the rectangle below. For $29$ and also for $31$ these multiples form a nice geometric picture.\label{29} }
\end{figure}
\begin{figure}[H]
\begin{center}
\includegraphics[scale=0.8]{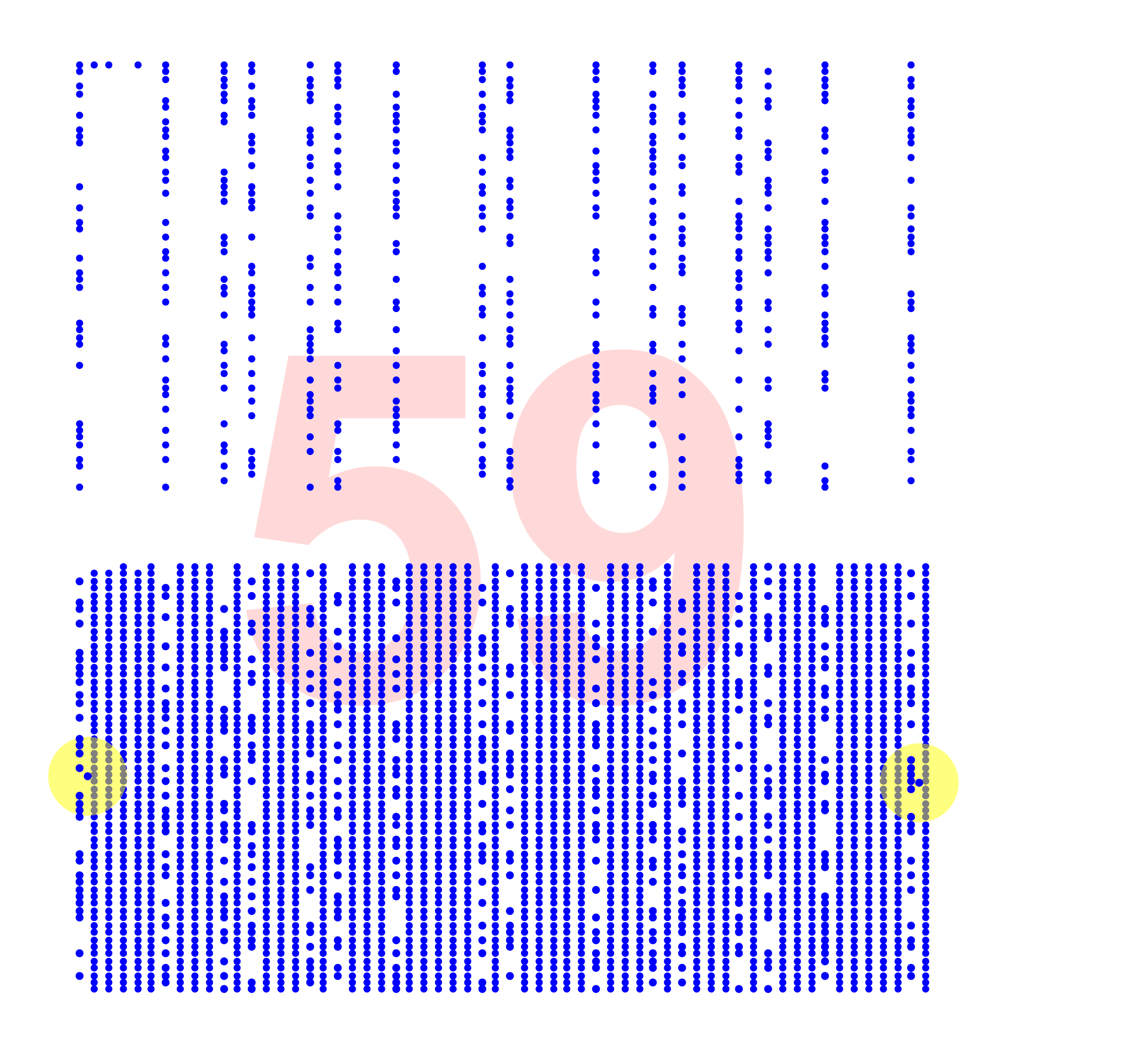}
\end{center}
\caption{All multiples of $59$ except $59$ itself are descending  in the rectangle below. The last two dancers correspond to $59^2$ on the left and $59\times 61$ on the right. \label{59} }
\end{figure}

\begin{figure}[H]
\begin{center}
\includegraphics[scale=0.8]{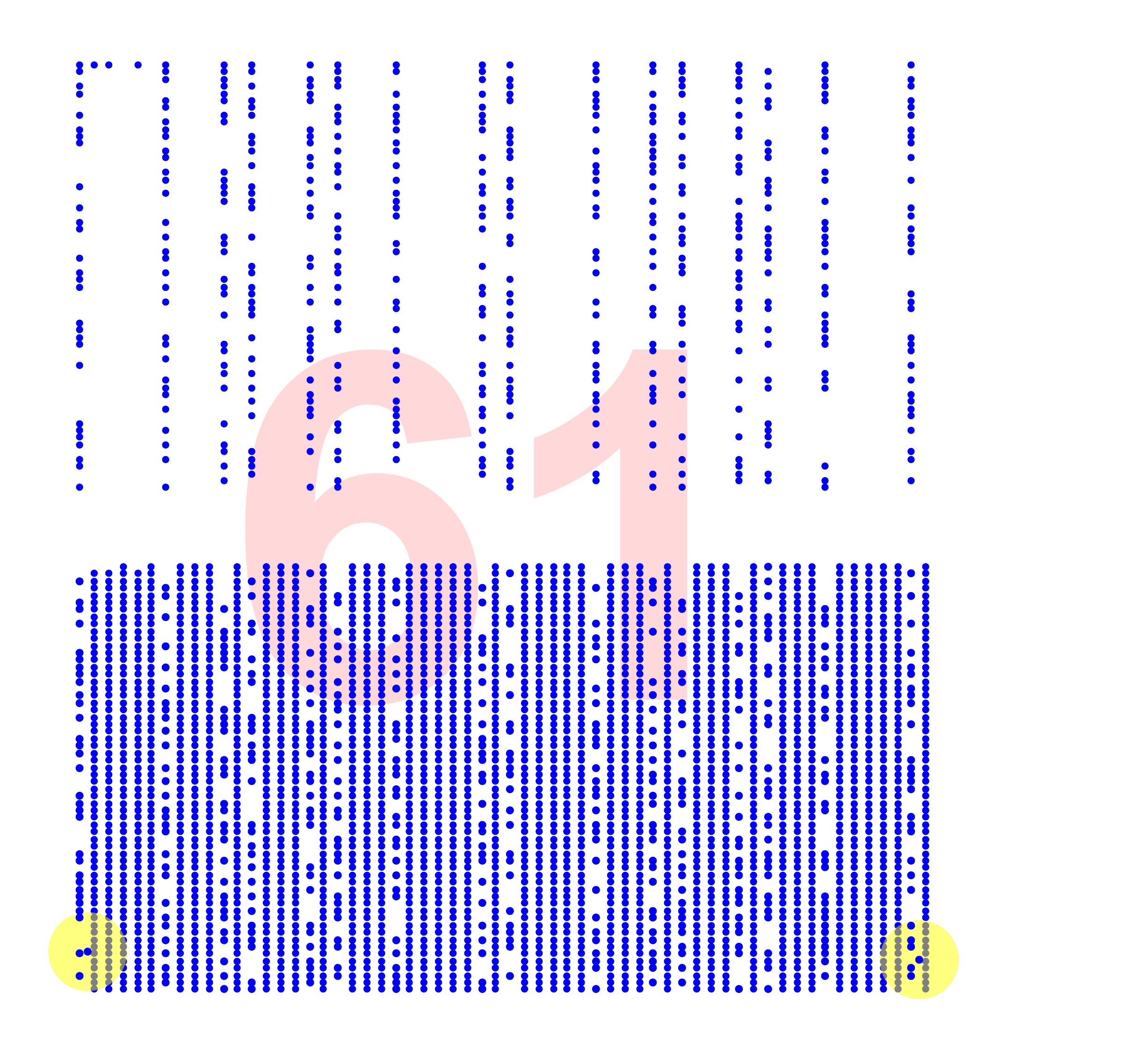}
\end{center}
\caption{ The last two dancers corresponding to $59^2$ on the left and $59\times 61$ on the right are about to reach their final position, thus ending the dance. \label{59} }
\end{figure}

\begin{figure}[H]
\begin{center}
\includegraphics[scale=1.3]{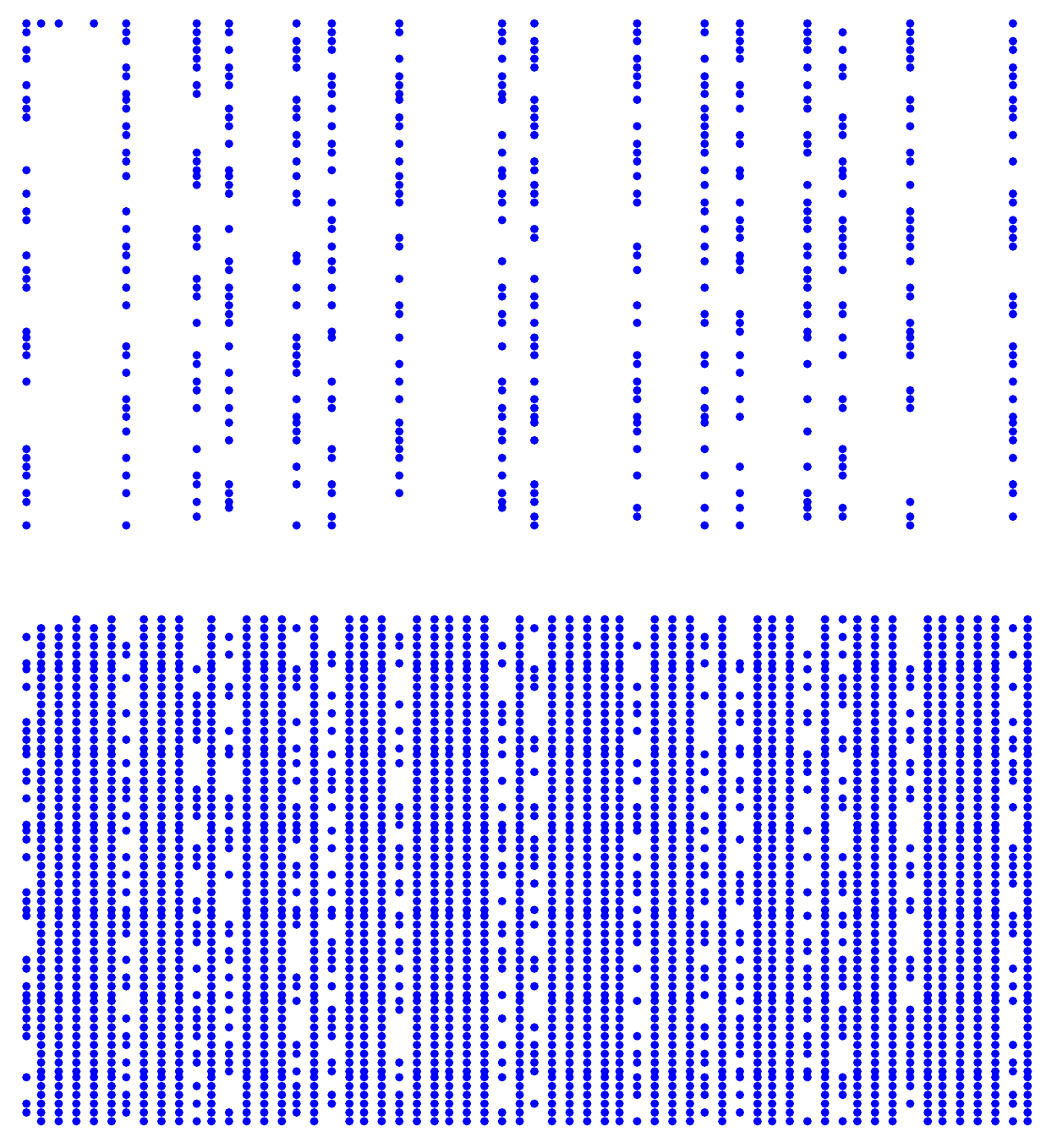}
\end{center}
\caption{After all multiples of $59$ except $59$ have descended, the dancers remaining in the top rectangle correspond exactly to the prime numbers below $3600$. What remains in the lower rectangle forms a sieve inasmuch as it is stable under multiplication by any number. In fact this property of being a sieve holds all along for the configurations in the lower rectangle.\label{end} }
\end{figure}

\subsection{Rhythmic interpretation by $H^1(C_j/p)$}

The second part of the video is an illustration of the rhythmic interpretation given by the eigenvalues of the Frobenius on $H^1(C_j/p)$, where $p$ is a prime $7\leq p \leq 67$
\begin{figure}[H]
\begin{center}
\includegraphics[scale=0.5]{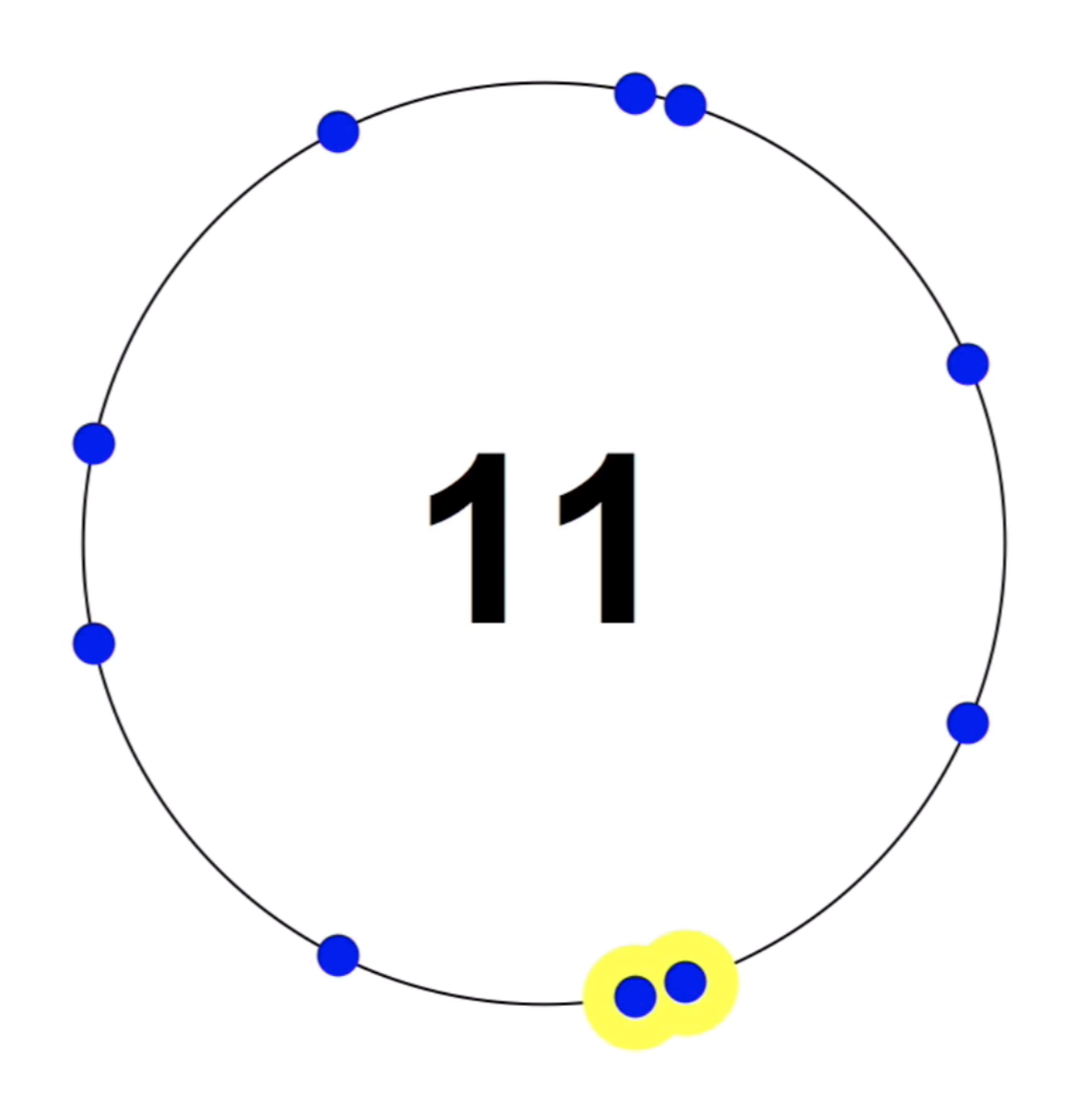}
\end{center}
\caption{Each of the six motives becomes an interpreter (at the rhythmic level) of the same piece.\label{rythm} }
\end{figure}
In order to illustrate these rhythms which provide the precise irrational times of attack in the music, the video is synchronized so that the corresponding zeros are shown in shining color while they are played. The tempo is accelerating with $\log p$ as $p$ varies from $7$ to $67$ and this is done for each of the six motives $H^1(C_j/p)$. The video is available by clicking on the following link:
 \href{https://www.dropbox.com/s/ld1z64em5u3mxjx/fullvideo.mp4?dl=0}{\color{blue}Video}

We have already described $C_1$ and $C_2$. The next ones and some samples of the associated rhythmic interpretations are as follows:

\newpage
$$
 C_3: \  y^2=x^{11}-x^{10}+x^9-5 x^8+8 x^7-8 x^6+8 x^5-14 x^4+5 x^3-7 x^2+x-1
 $$

 \begin{figure}[H]
\begin{center}
\includegraphics[scale=0.7]{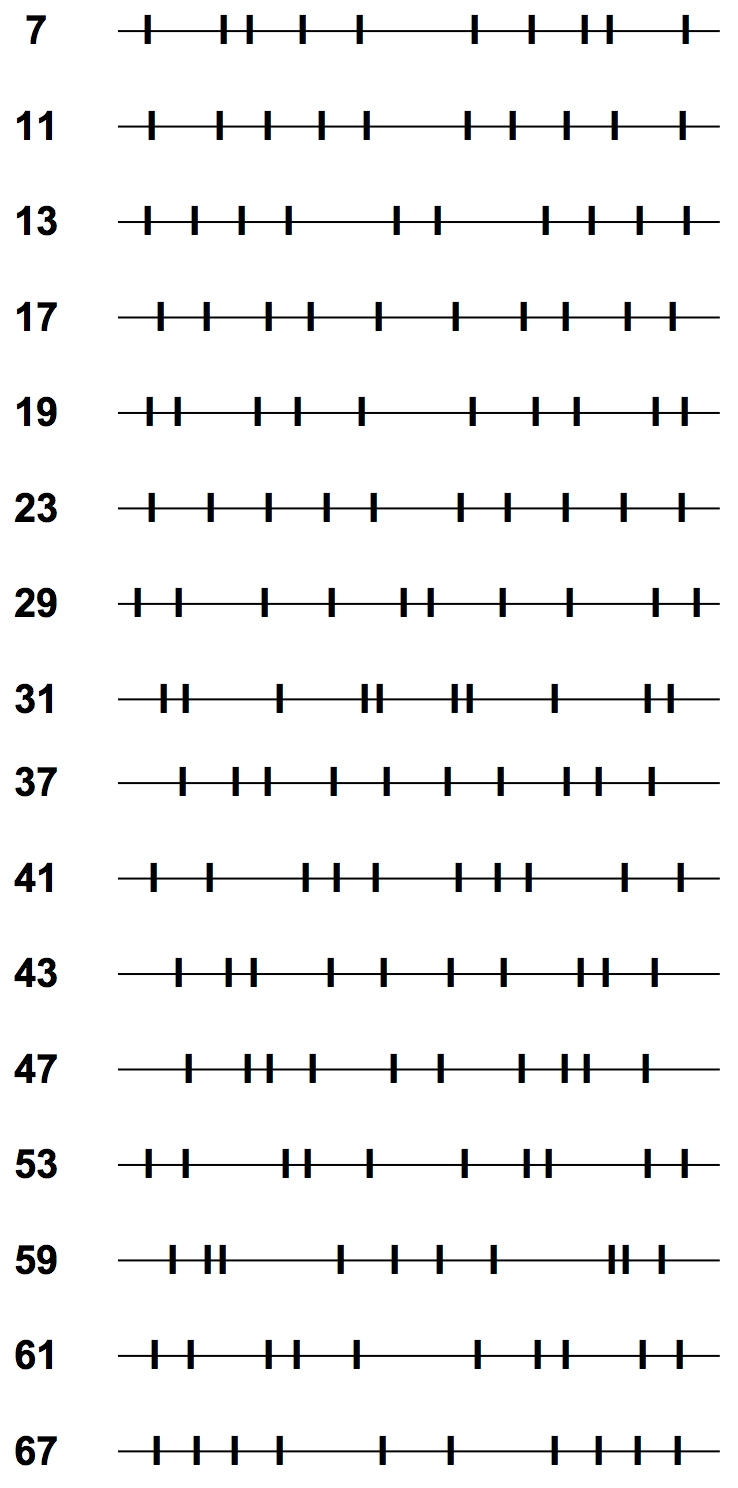}
\end{center}
\caption{Rhythms associated to $H^1(C_3/p)$ for $7\leq p\leq 67$. \label{poly1} }
\end{figure}
\newpage
 $$
C_4: \  y^2= x^{11}-x^{10}+3 x^9+x^8-8 x^7-8 x^5+24 x^4+58 x^3+86 x^2+86 x+50
 $$
\begin{figure}[H]
\begin{center}
\includegraphics[scale=0.7]{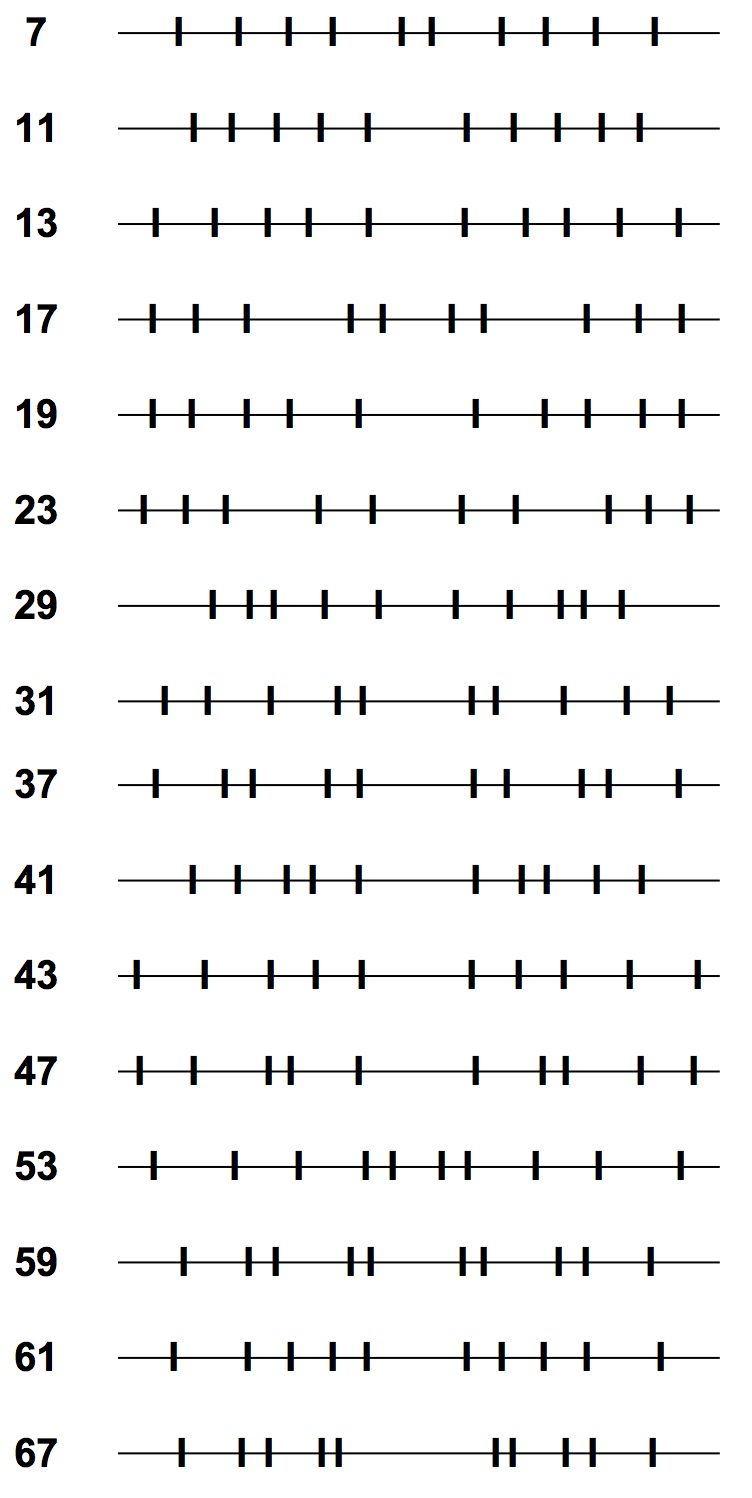}
\end{center}
\caption{Rhythms associated to $H^1(C_4/p)$ for $7\leq p\leq 67$. \label{poly1} }
\end{figure}
\newpage
$$
 C_5: \  y^2=x^{11}-x^{10}+7 x^9-15 x^8+36 x^7-48 x^6+108 x^5-144 x^4+90 x^3-162 x^2+162 x+198
 $$
\begin{figure}[H]
\begin{center}
\includegraphics[scale=0.7]{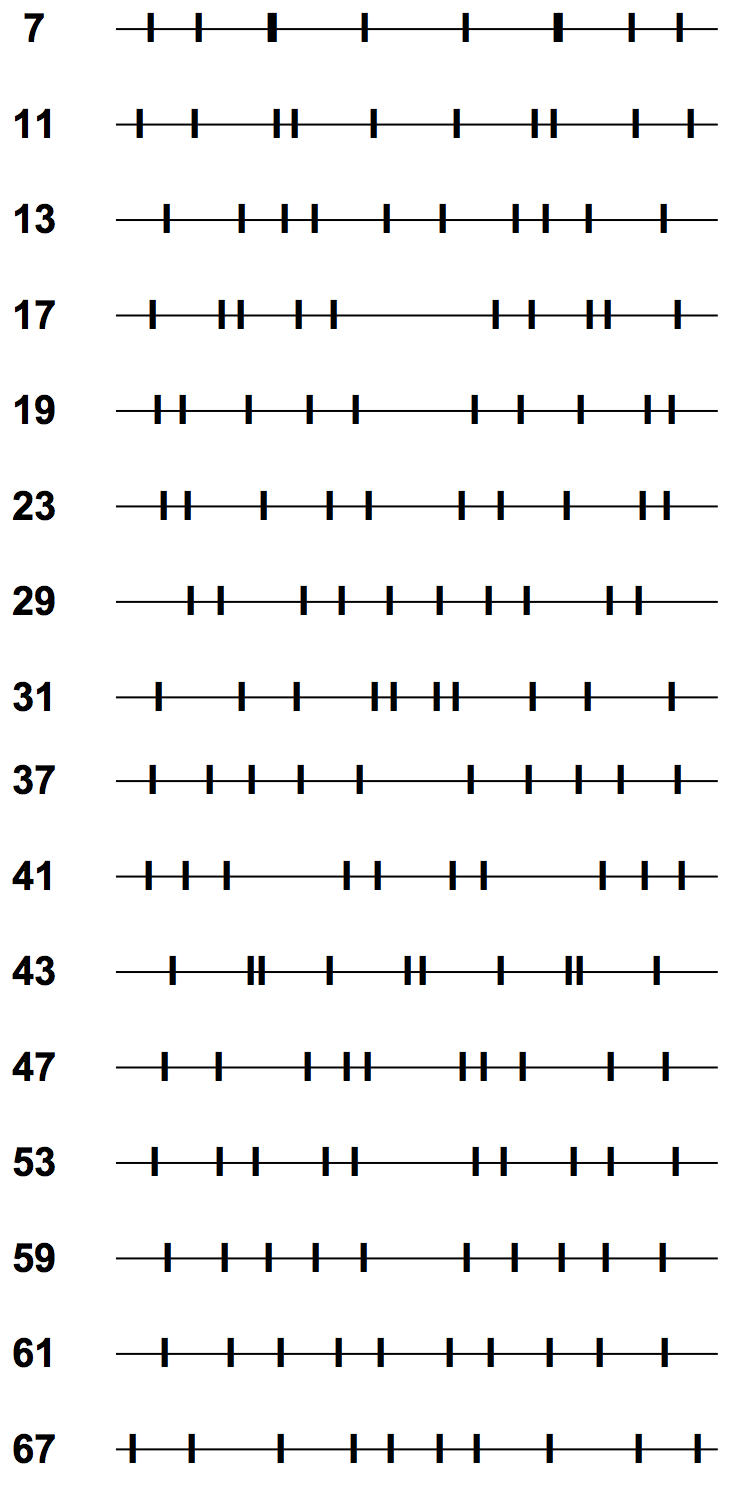}
\end{center}
\caption{Rhythms associated to $H^1(C_5/p)$ for $7\leq p\leq 67$. \label{poly1} }
\end{figure}
\newpage
$$
  C_6: \   y^2=x^{11}-3 x^{10}-6 x^9+12 x^8+18 x^7-54 x^6-96 x^5+72 x^4+126 x^3-206 x^2-336 x-96
 $$

\begin{figure}[H]
\begin{center}
\includegraphics[scale=0.7]{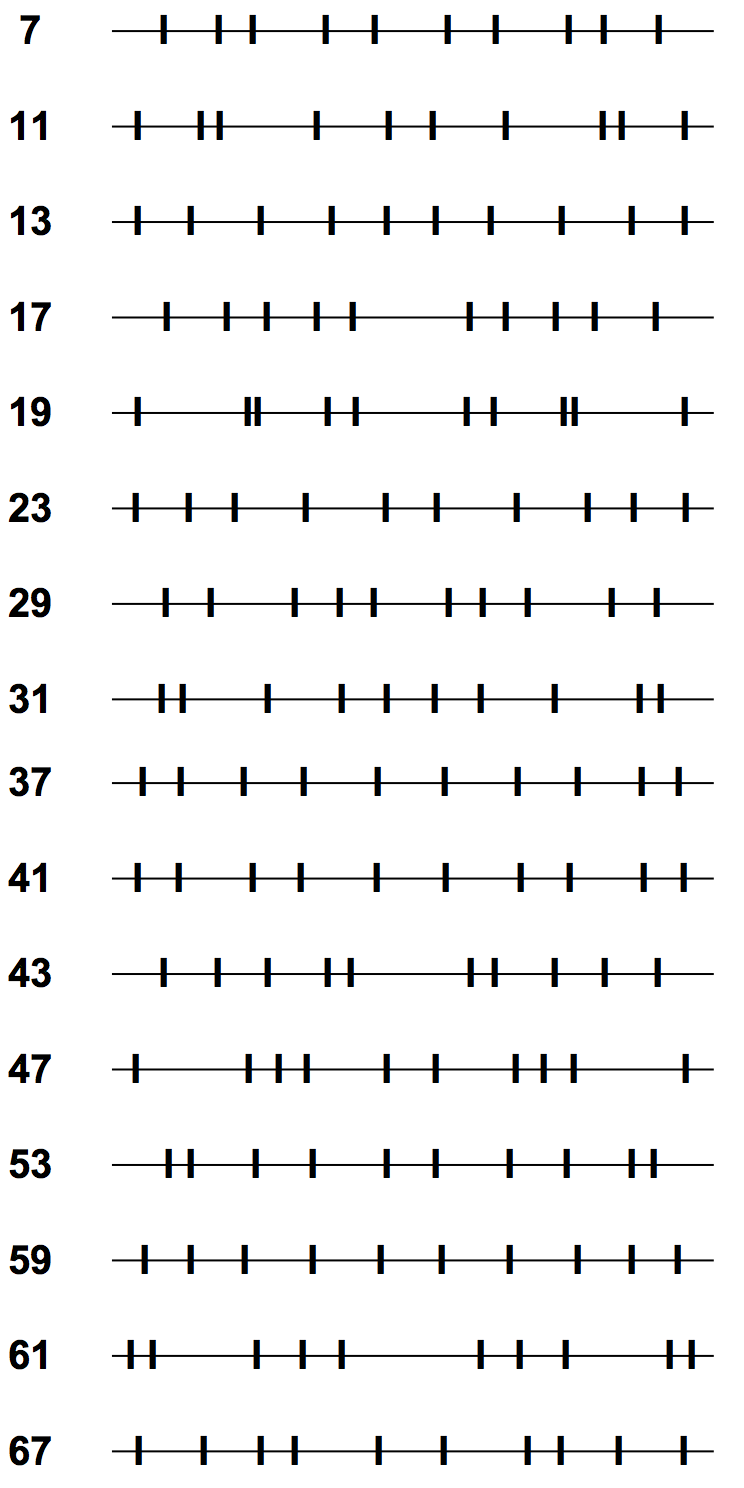}
\end{center}
\caption{Rhythms associated to $H^1(C_6/p)$ for $7\leq p\leq 67$. \label{poly1} }
\end{figure}
\newpage

\begin{rem} \label{finalrem} The fact that the tempo accelerates with $\log p$, \ie that the natural period determined by $p^{-is}=1$ is $\frac{2\pi}{\log p}$ exhibits in which sense the strategy followed in our joint work with C. Consani (see in particular \cite{CC0,CC0.5}) is natural. Indeed if one would let for instance $p\to \infty$ the period would tend to $0$ and there would be no way to get a picture resembling the zeros of the Riemann zeta function. But in our strategy one lets $p\to 1$ so that the period $\frac{2\pi}{\log p}$ now tends to infinity. Moreover as explained in \cite{CC0,CC0.5}, the genus tends to $\infty$, thus allowing for the desired type of configuration.
\end{rem}

\renewcommand{\abstractname}{Acknowledgements}
\begin{abstract}
 I am grateful to the referees for their useful comments, to Fernando Rodriguez Villegas,  Moreno Andreatta, Karim Haddad, Carlos Agon  and Ugo Moschella for their help in the realization of these ideas, and also to our IHES group of discussions on mathematics and music with, besides Moreno, Charles Alunni, Olivia Caramello and Pierre Cartier.  The principle of the above developments was obtained during the collaboration with Danye Ch\' ereau and Jacques Dixmier which appears in the novel \cite{CC1}.
\end{abstract}


\begin{thebibliography}{99}


\bibitem{Babbitt}  M.~Babbitt, \emph{ Twelve-Tone Rhythmic Structure and the Electronic Medium}, Perspectives of New Music 1, no. 1 (Fall 1962): 49--79

\bibitem{BKe} M.~Berry and J.~Keating, \emph{
 The Riemann zeros and eigenvalue asymptotics}. SIAM Rev. \textbf{41} (1999), no. 2, 236--266.

\bibitem{CC1} D.~Ch\' ereau, A.~Connes, J. Dixmier, {\em Le Spectre d'Atacama}. Editions Odile Jacob, 2018.



\bibitem{Co-zeta} A.~Connes, {\em Trace formula in noncommutative
geometry and the zeros of the Riemann zeta function}.  Selecta Math.
(N.S.)  5  (1999),  no. 1, 29--106.

\bibitem{CC0} A.~Connes, C.~Consani, {\em Schemes over $\F_1$ and zeta functions}, Compositio Mathematica
146 (6), (2010) 1383--1415.

\bibitem{CC0.5} A.~Connes, C.~Consani, {\em From monoids to hyperstructures: in search of an absolute arithmetic}, in Casimir Force, Casimir Operators and the Riemann Hypothesis, de Gruyter
(2010), 147--198.

\bibitem{Kac} M.~Kac, {\em Can One Hear the Shape of a Drum?}, The American Mathematical Monthly, Vol. 73, No. 4, Part 2: Papers in Analysis (Apr. 1966), pp. 1--23
(Available online at: https://www.math.ucdavis.edu/~hunter/m207b/kac.pdf)

\bibitem{Messiaen0} Olivier Messiaen,  {\em Technique de mon langage musical}, Paris, Alphonce Leduc, 1944 (English translation as "The Technique of My Musical Language", Paris, Alphonse Leduc, 1956.

\bibitem{Messiaen} O.~Messiaen, {\em Trait\' e de Rythme, de Couleur, et d'Ornithologie}. Editions musicales Alphonse Leduc.

\bibitem{Pressing} J.~Pressing, {\em Cognitive Isomorphisms between Pitch and Rhythm in World Musics: West Africa, the Balkans and Western Tonality}, Studies in Music, 17, (1983) 38--61.

\bibitem{Rhan} J.~Rahn, {\em On Pitch or Rhythm: Interpretations of Orderings of and in Pitch and Time,} Perspectives of New Music, vol. 13, no. 2, 1975, p. 182.
\end{thebibliography}
\end{document}